\documentclass[final]{siamart220329}
\usepackage{amsfonts, macrodefs, algorithmic, graphicx, subcaption, nicematrix}
\usepackage{xcolor}

\title{Structure-Aware Analyses and Algorithms for Interpolative Decompositions}

\author{Robin Armstrong\thanks{Cornell University, Center for Applied Mathematics, Ithaca, NY (\email{rja243@cornell.edu}).} \and Alex Buzali\thanks{Harvard University, School of Engineering and Applied Sciences, Cambridge, MA (\email{alexbuzali@g.harvard.edu}).} \and Anil Damle\thanks{Cornell University, Department of Computer Science, Ithaca, NY (\email{damle@cornell.edu}).}}

\begin{document}
\maketitle

\begin{abstract}
Low-rank approximation is a task of critical importance in modern science, engineering, and statistics. Many low-rank approximation algorithms, such as the randomized singular value decomposition (RSVD), project their input matrix into a subspace approximating the span of its leading singular vectors. Other algorithms compress their input into a small subset of representative rows or columns, leading to a so-called interpolative decomposition. This paper investigates how the accuracy of interpolative decompositions is affected by the structural properties of the input matrix being operated on, including how these properties affect the performance comparison between interpolative decompositions and RSVD. We give particular focus to a randomized Golub-Klema-Stewart (RGKS) algorithm, an interpolative  decomposition algorithm which combines RSVD with the Golub-Klema-Stewart (GKS) pivoting strategy. Through numerical experiments, we find that matrix structures including singular subspace geometry and singular spectrum decay play a significant role in determining the performance comparison between these different algorithms. We also prove inequalities which bound the error of a general interpolative decomposition in terms of these matrix structures. Finally, we develop forms of these bounds specialized to RGKS while considering how randomization affects the approximation error of this algorithm.
\end{abstract}



\section{Introduction}

A low-rank approximation of an $m \times n$ matrix $\mA$ is a decomposition
\begin{equation*}
\mA \approx \mB_1\mB_2\tp,
\end{equation*}
where $\mB_1$ is $m \times k$, $\mB_2$ is $n \times k$, and $k \ll \min\{ m,\, n \}$. We can refer to this more specifically as a rank-$k$ approximation. Approximations of this type are critical for efficiently handling extremely the extremely large matrices arising in modern applications. The key property of this approximation is that the dimensions of $\mB_1$ and $\mB_2$ are much smaller than those of $\mA$, meaning that by working with $\mB_1$ and $\mB_2$ instead of with $\mA$ directly, basic computations can be performed more efficiently.

\begin{figure}[t!]
    \centering
    \includegraphics[scale=0.7]{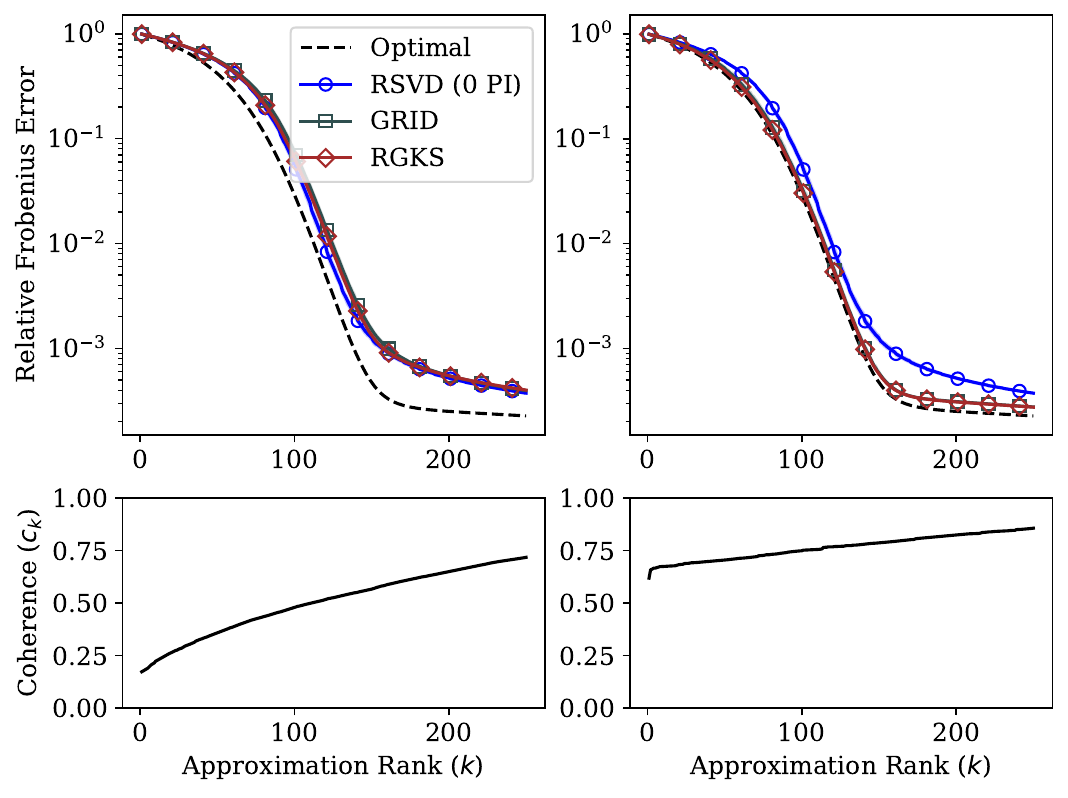}
    \caption{Relative Frobenius error for approximations of two $512 \times 512$ matrices computed by GRID, RGKS, RSVD with zero power iterations. Each curve shows the mean error over 100 approximations of the same matrix. The left-hand plot uses a test matrix with less coherent right singular subspaces, while the test matrix for the right-hand plot has more coherent subspaces.}
    \label{fig:algs_versusrank_intro}
\end{figure}

\begin{figure}[t!]
    \centering
    \includegraphics[scale=0.7]{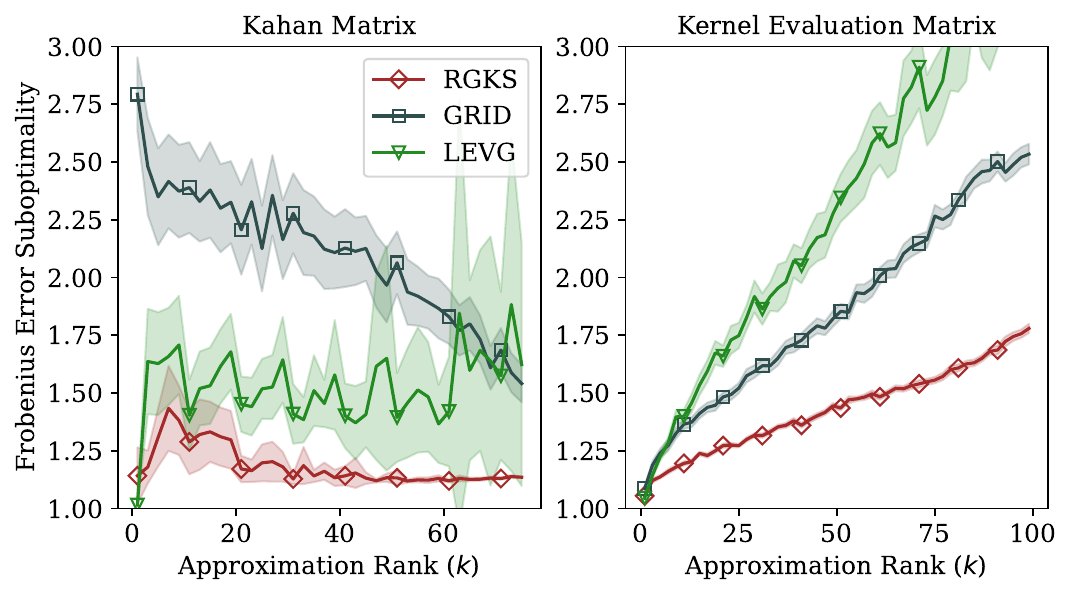}
    \caption{Frobenius error suboptimality factors for approximations computed by GRID, RGKS, and LEVG. Left panel: a $100 \times 100$ Kahan matrix. Right panel: an $800 \times 800$ matrix $\mA_{i,\, j} = \kappa(x_i,\, x_j)$, where $\kappa(x,\, y) = \exp(-\frac{1}{8}\| x - y \|_2^2)$, and $x_1,\, \ldots,\, x_{800}$ are i.i.d. samples from a 2D isotropic Gaussian with variance 25. Solid lines indicate a mean over 100 trials, and shaded regions indicate a 95\% confidence interval on the mean. Reasons for this behavior are discussed in \cref{section:numerics}.}
    \label{fig:kahan_kernel_example}
\end{figure}

Many algorithms for computing a rank-$k$ approximation fall into two broad categories. In the first category are algorithms that identify a ``structurally important'' $k$-dimensional subspace $\cQ \subseteq \R^m$, represented by an orthonormal basis $\mQ \in \R^{m \times k}$ whose columns are estimates of the leading left singular vectors of $\mA$--- leading to the low-rank approximation $\mB_1 = \mQ,\, \mB_2 = \mA\tp\mQ$. Algorithms in this category include randomized SVD, subspace iteration, and block-Krylov methods \cite{halko_finding_structure_with_randomness, musco_musco_rbki, rokhlin_randomized_pca, tropp_2023_randomized}. We will focus primarily on the randomized SVD (RSVD). The second category encompasses so-called interpolative and CUR decompositions, which identify a small subset of ``skeleton'' indices $J = \{ j_1,\, \ldots,\, j_k \}$ corresponding to the most ``structurally important'' rows and/or columns. In the case where $J$ is a set of column indices, this leads to a low-rank approximation $\mB_1 = \mA_{:,\, J},\, \mB_2\tp = (\mA_{:,\, J})\pinv\mA$, where the subscript ``$:,\, J$'' is \textsc{Matlab} notation for selecting a column subset, and $\dagger$ denotes the Moore-Penrose pseudoinverse. This paper will focus on approximations based on column selection. Minimizing the approximation error over all $n!/(k!(n - k)!)$ choices of $J$ is NP hard \cite{civril_2009}, but a variety of algorithms exist for choosing approximately optimal columns. These include approaches based on pivoted QR or LU factorizations applied to the columns of $\mA$ \cite{chandrasekaran_1994, cheng_interpolative_decomposition, dong_simpler_is_better, gu_srrqr, hong_pan, martinsson_rokhlin_tygert_randid, woolfe_srft}, random column sampling strategies \cite{deshpande_rademacher_volumesampling, deshpande_rademacher_projectiveclustering, mahoney_drineas_muthukrishnan, frieze_monte_carlo, mahoney_cur}, and selection strategies based on the singular value decomposition \cite{mahoney_drineas_muthukrishnan, drmac_qdeim, golub_klema_stewart, mahoney_cur, sorensen_embree_deim}.
Although most matrices do not have a $k$-column subset that exactly spans the optimal (i.e.\ leading singular) subspace, there may nevertheless exist a column subset whose span closely approximates this subspace. In these situations, using randomization to identify these columns may produce a more accurate approximation than directly estimating the leading singular subspace by sketching. To convey this point, we refer to \cref{fig:algs_versusrank_intro}, which plots the relative approximation error of three randomized low-rank approximation algorithms as a function of the approximation rank. One of these algorithms is the RSVD, and the other two are interpolative decompositions: GRID (\cref{alg:grid}), which uses the well-known technique of applying a pivoted QR factorization to a randomized embedding of the input columns \cite{deutsch_gu_rqrcp, martinsson_rokhlin_tygert_randid, martinsson_hqrrp, woolfe_srft}, and RGKS (\cref{alg:rgks}), a randomized variant of the Golub-Klema-Stewart algorithm \cite{golub_klema_stewart}. An important fact illustrated by \cref{fig:algs_versusrank_intro} is that for certain test matrices, the GRID and RGKS have accuracy which is competitive with RSVD when power iteration is not used, a comparison which places all algorithms on a similar footing in terms of the number of matrix-vector products with $\mA$. And although GRID and RGKS perform essentially identically in \cref{fig:algs_versusrank_intro}, the structure of a matrix can significantly change the performance comparison between different column selection strategies. For example, \cref{fig:kahan_kernel_example} shows the approximation error of GRID, RGKS, and leverage score sampling (\cref{alg:levg}) for a Kahan matrix,\footnote{While the $n \times n$ Kahan matrix is typically considered a pathological test case for column-pivoted QR factorizations, the well-known pathological behavior occurs at $k = n - 1$. The behavior in \cref{fig:kahan_kernel_example} is more surprising as we are considering smaller values of $k$, and we would expect the randomness to partially mitigate the delicate structure that ``breaks'' the standard column-pivoted QR algorithm.} as well as a matrix generated by evaluating a Gaussian kernel pairwise between points randomly drawn from a two-dimensional Gaussian distribution. For these matrices, in contrast to the test matrices of \cref{fig:algs_versusrank_intro}, RGKS performs significantly better than its counterparts. Reasons for this are discussed in \cref{section:numerics}.

Beyond simply highlighting that interpolative decompositions can be competitive with RSVD given a fixed budget of matrix-vector products, and that different column selection strategies are more well suited to different matrix structures, this paper seeks to characterize specific structures of the input matrix which affect the comparative accuracy of these algorithms. The structures affecting approximation accuracy which we examine can be categorized into two groups: those which characterize the decay of a matrix's singular spectrum, and those which characterize the geometry of its singular subspaces. For spectral decay, we will focus on the singular value gap and residual stable rank, defined by
\begin{equation*}
\gamma_k = \frac{\sigma_{k + 1}}{\sigma_k} \quad\text{and}\quad r_k = \sum_{i = k + 1}^{\min\{ m,\, n \}} \frac{\sigma_i^2}{\sigma_{k + 1}^2},
\end{equation*}
where $\sigma_1 \geq \ldots \geq \sigma_{\min\{ m,\, n \}} \geq 0$ are the singular values of $\mA$.

To describe the geometry of a matrix's singular subspaces, we will use several concepts related to the singular value decomposition $\mA = \mU\mSigma\mV\tp$. The central quantity of interest in this context is $\sigma_\mathrm{min}(\mV_{J,\, 1:k})$, where $J$ is the set of skeleton column indices. Many results connecting $\sigma_\mathrm{min}(\mV_{J,\, 1:k})$ to interpolative decomposition accuracy have appeared in the literature before; for example, it is known that
\begin{equation}
\| \mA - \mA_{:,\, J}(\mA_{:,\, J})\pinv\mA \|_2 \leq \frac{\sigma_{k + 1}(\mA)}{\sigma_\mathrm{min}(\mV_{J,\, 1:k})} \quad\text{and}\quad \sin \phi_\mathrm{max} \leq \frac{\gamma_k}{\sigma_\mathrm{min}(\mV_{J,\, 1:k})}, \label{eqn:prior_sigmamin_bounds}
\end{equation}
where $\phi_\mathrm{max}$ is the largest principal angle between $\range(\mA_{:,\, J})$ and $\range(\mU_{:,\, 1:k})$. These inequalities arise from \cite[Theorem 1.5]{hong_pan} and \cite[Theorems 3.1, 6.1]{golub_klema_stewart}, respectively. This paper examines how $\sigma_\mathrm{min}(\mV_{J,\, 1:k})$ is related to geometric properties of the right singular subspace, namely, the distribution of the subspace's column-based leverage scores. The graphs in \cref{fig:algs_versusrank_intro} differ in terms of the leverage score distribution of the input matrix's right singular subspaces, and from this figure, we see that this difference strongly affects the performance of the two interpolative decompositions. We will also find that assigning a geometric interpretation to $\sigma_\mathrm{min}(\mV_{J,\, 1:k})$ allows several error bounds to be presented in a unified notation, and simplifies our analysis of randomization errors in the RGKS algorithm (described below).

Related to our structural analysis of interpolative decompositions, this paper gives special focus to the \emph{randomized Golub-Klema-Stewart algorithm }(RGKS), a randomized variant of GKS \cite{golub_klema_stewart}, which is conceptually related to the QDEIM algorithm that later appeared in model order reduction \cite{drmac_qdeim}. At a high level, RGKS selects skeleton columns by applying a rank-revealing QR factorization to a matrix of right singular vector estimates computed by the RSVD.  Because RGKS operates directly on an estimate of the right singular subspace, it serves as a natural case study on the effects of leverage score distributions. Furthermore, because the subspace estimate is computed using RSVD, whose accuracy is known to depend on gaps in the singular spectrum, it will be natural to consider how the spectral and subspace effects couple with one another. In particular, prior work analyzes RSVD subspace estimates, e.g., \cite{dong_martinsson_nakatsukasa, saibaba_randomized_subspace_iteration, zhang_tang}, and this raises interesting possibilities for exploring how the particular subspace errors committed by RSVD interact with the subspace-dependent QR factorization in RGKS. Saibaba's work on the RDEIM algorithm \cite{saibaba_rdeim} is related to our work on RGKS, and relevant differences between these two algorithms and associated theoretical results are explained in \cref{section:rgks}.

While RGKS serves as an ``archetypal'' interpolative decomposition to motivate our analysis of interpolative decompositions in general, we will also consider several interesting features which are unique to RGKS. Specifically, we will show that RGKS has attractive properties in terms of efficiency, accuracy, and robustness to noise arising from its internal randomization. We will also show that in many cases, interpolative decompositions computed with RGKS are more accurate than those computed with leverage score sampling, a randomized column subset selection procedure which also selects skeleton indices based on the right singular subspace (or an approximation thereof).

\subsection{Notes on Numerical Experiments}

Except where otherwise stated, all instances of RSVD in this paper were run without power iteration. All low-rank approximation algorithms were run with oversampling equal to $p = \lceil k/10 \rceil$, where $k$ is the target approximation rank. Our experiments will mostly measure errors in the Frobenius norm; supplementary experiments measuring spectral norm error are included in the supplementary text.

\subsection{Main Contributions}

The main algorithmic contribution of this paper is RGKS, a sketching-based interpolative decomposition which combines a randomized SVD with a rank-revealing QR factorization. \Cref{section:rgks} describes this algorithm, and discusses its accuracy and efficiency in comparison to other interpolative decomposition algorithms. \Cref{section:rgks_analysis} develops error bounds for RGKS, and shows that the accuracy of RGKS is surprisingly robust to singular vector estimation errors arising from randomization.

The analytical contributions of this paper include error bounds that characterize how the structural properties of an input matrix affect the accuracy of interpolative decompositions. In stating these bounds, $\mSigma_\perp = \mathrm{diag}(\sigma_i \suchthat i \geq k + 1)$ denotes the matrix of trailing singular values of $\mA$, so that $\| \mSigma_\perp \|$ is the optimal approximation error in the spectral or Frobenius norm. The span of the first $k$ right singular vectors of $\mA$ is denoted by $\cV_k$, and to ensure that this subspace is well defined, we assume that $\sigma_k(\mA) > \sigma_{k + 1}(\mA)$. As is customary, $\ve_j$ denotes the $j\nth$ elementary unit vector. We now summarize our main error bounds in the theorem below, which is a concatenation of \cref{thm:idbound_subspaceonly,thm:idbound_bothstructures,thm:lowrankbound_spectrumonly}.

\begin{theorem}
Choose $k \leq n/2$ such that $\sigma_k(\mA) > \sigma_{k + 1}(\mA)$. Given a set of skeleton column indices $J \subseteq \{ 1,\, \ldots,\, n \}$ with $|J| = k$, define the approximation error $\mE = \mA - \mA_{:,\, J}(\mA_{:,\, J})\pinv\mA$, and let $\varphi_1,\, \ldots,\, \varphi_k \in [0,\, \pi / 2]$ be the principal angles between $\cV_k$ and $\cI_J \defeq \vspan\{ \ve_j \suchthat j \in J \}$. If $\max_{1 \leq i \leq k} \varphi_i < \pi / 2$, then
\begin{equation}
\| \mE \|_2 \leq \| \mSigma_\perp \|_2 \sec ({\textstyle \max_{1 \leq i \leq k} \varphi_i}) \quad \text{and} \quad \| \mE \|_\frob \leq \| \mSigma_\perp \|_\frob \sqrt{1 + \frac{1}{r_k} \sum_{i = 1}^k \tan^2 \varphi_i}, \label{eqn:subspacebounds_intro}
\end{equation}
where $r_k = \| \mSigma_\perp \|_\frob^2 / \| \mSigma_\perp \|_2^2$ is the residual stable rank. If $\sigma_{k + 1}(\mE) > 0$ then, in addition,
\begin{equation}
\| \mE \|_2 \leq \| \mSigma_\perp \|_2 \kappa(\mE,\, k + 1), \label{eqn:conditionbound_intro}
\end{equation}
where $\kappa(\mE,\, k + 1) = \sigma_1(\mE)/\sigma_{k + 1}(\mE)$ is a modified condition number.
\end{theorem}

In \cref{section:analysis}, we will show that the angles $\varphi_1,\, \ldots,\, \varphi_k$ encode the conditioning of the row-subset $\mV_{J,\, 1:k}$, and more generally the effects of subspace geometry on interpolative decomposition accuracy. The stable rank $r_k$ in \cref{eqn:subspacebounds_intro} encodes the effects of singular spectrum decay for Frobenius norm errors, and we will show that when the trailing singular value spectrum is nearly flat, \cref{eqn:conditionbound_intro} provides an exceptionally tight bound on spectral norm errors. For reasons explained in \cref{subsec:idbound_subspaceonly}, the spectral norm bound in equation \cref{eqn:subspacebounds_intro} is equivalent to a 1992 result of Hong and Pan \cite[Theorem 1.5]{hong_pan}. Our formulation of the result emphasizes a geometric interpretation and a connection to interpolative decompositions, whereas Hong and Pan stated the result as a singular value inequality for rank-revealing QR factorizations. Sorensen and Embree stated and generalized the same singular value inequality in the context of matrix approximation \cite{sorensen_embree_deim}. Our formulation of the result in terms of principal angles allows us to present several different interpolative decomposition error bounds in a unified notation, and also simplifies our analysis of randomization errors in RGKS.

Finally, this paper provides a set of numerical experiments which compare the accuracy of RGKS, leverage score sampling, and RSVD across variations of structure in the input matrix, mainly in \cref{section:numerics}. These experiments demonstrate the potential of interpolative decompositions to have competitive or superior accuracy to RSVD, given certain choices of approximation rank and certain structural properties of the input matrix. This goes against the common intuition that the RSVD, being an approximation to the optimal SVD, should always have superior performance relative to the more ``restrictive'' interpolative decompositions.

\section{Background and Notation} \label{section:background}

Computing a low-rank approximation amounts to solving, either exactly or approximately, the minimization problem
\begin{equation}
\min \{ \| \mA - \mB_1 \mB_2\tp \| \suchthat \mB_1 \in \R^{m \times k},\, \mB_2 \in \R^{n \times k} \}, \label{eqn:lowrank_minproblem}
\end{equation}
where $\| \cdot \|$ is a matrix norm. Many aspects of this problem are determined by the singular value decomposition $\mA = \mU\mSigma\mV\tp$, especially when it is partitioned at rank $k$ in the following manner:
\begin{equation}
\mA = \begin{bmatrix} \mU_k & \mU_\perp \end{bmatrix} \begin{bmatrix}
\mSigma_k & \mZero \\
\mZero & \mSigma_\perp
\end{bmatrix} \begin{bmatrix}
\mV_k & \mV_\perp
\end{bmatrix}\tp = \mU_k\mSigma_k\mV_k\tp + \mU_\perp\mSigma_\perp\mV_\perp\tp,
\end{equation}
where $\mSigma_k = \mathrm{diag}(\sigma_1,\, \ldots,\, \sigma_k)$ contains the largest $k$ singular values of $\mA$, and $\mSigma_\perp = \mathrm{diag}(\sigma_{k + 1},\, \ldots,\, \sigma_{\min\{ m,\, n \}})$ contains the remaining ones. The subspaces spanned by the leading singular vectors are
\begin{equation*}
\cU_k = \range(\mU_k),\qquad \cV_k = \range(\mV_k),
\end{equation*}
and we will refer to these as the leading singular subspaces. The significance of the SVD in low-rank approximation comes from the Eckart-Young-Mirsky theorem \cite{eckart_young_lowrank_optimal, mirsky_lowrank_optimal}, which states that if $\| \cdot \|$ is any unitarily invariant norm, then an optimal solution to \cref{eqn:lowrank_minproblem} can be obtained by projecting $\mA$ column-wise into $\cU_k$ (i.e., setting $\mB_1 = \mU_k,\, \mB_2 = \mA\tp\mU_k$), or by projecting row-wise into $\cV_k$ ($\mB_1 = \mA\mV_k,\, \mB_2 = \mV_k$).

In light of this result, many low-rank approximation algorithms attempt to form a subspace approximating $\cU_k$. Various algorithms in this category proceed by taking matrix-vector products between $\mA$ and a small number of randomized  test vectors, producing a collection of output vectors  termed a ``sketch'' of $\mA$. This is used to form an approximation
\begin{equation*}
\mA \approx \widehat{\mU}_k\widehat{\mSigma}_k\widehat{\mV}_k,\quad \widehat{\mSigma}_k = \mathrm{diag}(\widehat{\sigma}_1,\, \ldots,\, \widehat{\sigma}_k),
\end{equation*}
where $\widehat{\sigma}_1 \geq \ldots \geq \widehat{\sigma}_k \geq 0$ are approximate singular values, and $\widehat{\mU}_k \in \R^{m \times k},\, \widehat{\mV}_k \in \R^{n \times k}$ are orthonormal bases for subspaces approximating $\cU_k$ and $\cV_k$. For these algorithms to accurately approximate $\cU_k$ and $\cV_k$, it is favorable for there to be a large gap between $\sigma_k(\mA)$ and $\sigma_{k + 1}(\mA)$ \cite{halko_finding_structure_with_randomness, saibaba_randomized_subspace_iteration}. Randomized subspace iteration algorithms \cite{halko_finding_structure_with_randomness, rokhlin_randomized_pca} amplifies this gap by implicitly forming a sketch of $(\mA\mA\tp)^q\mA$, where $q \geq 0$ is a ``power iteration'' parameter that leaves the singular subspaces invariant while exponentially increasing the rate of singular value decay. Randomized block-Krylov algorithms \cite{musco_musco_rbki, rokhlin_randomized_pca, tropp_2023_randomized} implicitly form a sketch of $p_q(\mA)$, where $p_q$ is a degree-$q$ polynomial filter that separates $\sigma_k(\mA)$ from $\sigma_{k + 1}(\mA)$. In this paper we focus on randomized subspace iteration, which we refer to by its more common name, the randomized singular value decomposition (RSVD). Pseudocode for RSVD is given in \cref{alg:rsvd}, which also inclues an ``overampling'' parameter $p$. In our implementation of RSVD, an oversampled factorization is always optimally reduced to rank $k$.
\begin{algorithm}
\caption{\textsc{Randomized Singular Value Decomposition (RSVD)}} \label{alg:rsvd}
\begin{algorithmic}[1]
\STATE $\mOmega \leftarrow$ an $n \times (k + p)$ matrix with i.i.d. standard Gaussian entries.
\STATE $\mQ,\, \mR \leftarrow \texttt{thin\_qr}(\mA\mOmega)$.
\FOR{$i = 1,\, \ldots,\, q$}
\STATE $\mQ,\, \mR \leftarrow \texttt{thin\_qr}(\mA\tp\mQ)$.
\STATE $\mQ,\, \mR \leftarrow \texttt{thin\_qr}(\mA\mQ)$
\ENDFOR
\STATE $\mU,\, \mSigma,\, \mV \leftarrow \texttt{thin\_svd}(\mQ\tp\mA)$.
\RETURN $\widehat{\mU}_k = \mQ\mU_{:,\,1:k},\, \widehat{\mSigma}_k = \mSigma_{1:k,\, 1:k},\, \widehat{\mV}_k = \mV_{:,\, 1:k}$.
\end{algorithmic}
\end{algorithm}

Interpolative decompositions are a distinct class of low-rank approximations which are formed in terms of row and/or column submatrices of $\mA$. Approximations of this sort inherit more structure from $\mA$ and tend to be more interpretable \cite{mahoney_cur, voronin_martinsson_efficient_cur}. Most basic is the column-based interpolative decomposition,
\begin{equation}
\mA \approx \mA_{:,\, J}\mX, \label{eqn:column_based_id}
\end{equation}
where $J$ is a subset of $k$ column indices, and $\mX \in \R^{k \times n}$ is a matrix such that $\mX_{:,J}$ is a column subset of the identity matrix. This approximation is ``interpolating'' in the sense that $(\mA_{:,\, J}\mX)_{:,\, J} = \mA_{:,\, J}$ exactly.

Two-sided interpolative decomposition also exist with the form $\mA \approx \mW\mA_{I,\, J}\mX$, where $I$ contains $k$ row indices and $\mW \in\ \R^{m \times k}$. Such an approximation can be built from two column-based decompositions, namely, $\mA \approx \mA_{:,\, J}\mX$ followed by $(\mA_{:,\, J})\tp \approx (\mA_{:,\, J})\tp_{:,\, I}\mW\tp$ \cite{cheng_interpolative_decomposition, voronin_martinsson_efficient_cur}. Alternately, if $\mA \approx \mA_{:,\, J}\mX$ and $\mA\tp \approx (\mA\tp)_{:,\, I}\mW\tp$ are computed separately, then with various straightforward methods \cite{mahoney_cur, voronin_martinsson_efficient_cur}, one can form a two-sided interpolating matrix $\mU \in \R^{k \times k}$ such that $\mA \approx \mA_{:,\, J}\mU\mA_{I,\, :}$, a so-called ``CUR'' decomposition. In all of these methods, the main building block is the column-based interpolative decomposition, and henceforth, we will only consider interpolative decompositions of this sort.

To compute a column-based interpolative decomposition, the main problem is to select the column indices $J$ so as to minimize approximation error. This problem is NP-hard \cite{civril_2009}, but various deterministic and randomized algorithms exist which choose a nearly-optimal $J$, based on the principle that the columns of $\mA_{:,J}$ should be large and well conditioned. Most deterministic algorithms are based on pivoted matrix factorizations, and operate either on the columns of $\mA$ or on a leading subset of its singular vectors. A widely used strategy is to apply a column-pivoted QR (CPQR) factorization to $\mA$'s columns, producing a decomposition 
\begin{equation}
\mA\mPi = \mQ \mR = \mQ \begin{bNiceMatrix}[last-row, last-col]
\mR_{11} & \mR_{12} & \mbox{\scriptsize $k$} \\
\mZero & \mR_{22} & \mbox{\scriptsize $\min\{m,\, n\} - k$} \\
\mbox{\scriptsize $k$} & \mbox{\scriptsize $n - k$}
\end{bNiceMatrix}\:\:, \label{eqn:cpqr}
\end{equation}
where $\mPi$ is a permutation matrix, $\mQ$ has orthonormal columns, and $\mR_{11}$ is upper-triangular. Taking $\mA_{:,\,J}=\mA\mPi_{:,\, 1:k}$ and $\mX = \begin{bmatrix} \mI & \mR_{11}^{-1}\mR_{12} \end{bmatrix}\mPi\tp$ results in an interpolative decomposition that satisfies
\begin{equation}
\| \mA - \mA_{:,J}\mX \| = \| \mR_{22} \|,\label{eqn:general_cpqr_id}
\end{equation}
where $\| \cdot \|$ is the spectral or Frobenius norm. Minimizing $\| \mR_{22} \|$ can be accomplished using rank-revealing QR (RRQR) factorization algorithms \cite{chandrasekaran_1994, gu_srrqr, hong_pan}, which are algorithms that compute \cref{eqn:cpqr} in such a way that
\begin{equation}
\sigma_\mathrm{min}(\mR_{11}) \geq \frac{\sigma_k(\mA)}{q(n,\, k)} \quad \text{and} \quad \sigma_\mathrm{max}(\mR_{22}) \leq q(n,\, k) \sigma_{k + 1}(\mA), \label{eqn:rrqr_bounds}
\end{equation}
where $q$ is a function whose growth in $n$ is bounded by a low-degree polynomial. The exact form of $q$ depends on which RRQR algorithm is used to compute \cref{eqn:cpqr}.

Alternately, in the model reduction community the DEIM algorithm \cite{chaturantabut_sorensen_deim} chooses $J$ using a row-pivoted LU-like procedure applied to $\mV_k$, while the later-developed QDEIM \cite{drmac_qdeim} applies a row-pivoted LQ factorization to $\mV_k$. DEIM can be used to compute interpolative decompositions where $\mX$ corresponds to an oblique projection of $\mA$'s columns into $\range(\mA_{:,J})$ \cite{sorensen_embree_deim}. An interpolative decomposition using the QDEIM selection strategy first appeared in a technical report by Golub, Klema, and Stewart \cite{golub_klema_stewart} several decades before the invention of QDEIM itself; this algorithm is therefore referred to as GKS. These algorithms are all motivated by results \cite{golub_klema_stewart, hong_pan, sorensen_embree_deim} which bound the error of \cref{eqn:column_based_id} in terms of $\sigma_\mathrm{min}(\mV_{J,\, 1:k})$.

Pivoting-based interpolative decompositions can be accelerated using randomization in various ways. A particularly well-known method, sometimes called ``RID,'' involves applying a CPQR factorization to $\mOmega\tp\mA$, where $\mOmega \in \R^{m \times k}$ is a random sketching matrix defining an approximate isometry from $\R^m$ to $\R^k$. The original descriptions of RID \cite{martinsson_rokhlin_tygert_randid, woolfe_srft} used a subsampled randomized Fourier transform for $\mOmega$, allowing for a very rapid sketch, and computed the interpolating matrix $\mX$ based on the embedded columns. The experiments in this paper will use a different implementation, wherein a Gaussian random matrix is used for $\mOmega$, and $\mX$ is recomputed based on the columns of $\mA$ itself. Our implementation of this algorithm, which we refer to as ``Gaussian RID'' (GRID), is given in \cref{alg:grid}, where $\texttt{rrqr}(\cdot,\, k)$ indicates a rank-revealing QR factorization for rank $k$. GRID is related to algorithms that use sketching to efficiently compute a CPQR factorization by blocking the permutation computation and using rapid approximate column-norm updates \cite{deutsch_gu_rqrcp} and/or BLAS-3 Householder reflections \cite{martinsson_hqrrp}.
\begin{algorithm}
\caption{\textsc{Gaussian RID (GRID)}} \label{alg:grid}
\begin{algorithmic}[1]
\STATE $\mOmega \leftarrow$ an $m \times (k + p)$ matrix with i.i.d. standard Gaussian entries.
\STATE $\mS \leftarrow \mOmega\tp\mA$.
\STATE $\mPi,\, \mQ,\, \mR \leftarrow \texttt{rrqr}(\mS,\, k)$.
\RETURN $\mB_1 = \mA\mPi_{:,\, 1:k},\, \mB_2\tp = (\mA\mPi_{:,\, 1:k})\pinv\mA$.
\end{algorithmic}
\end{algorithm}

Another widely studied technique for skeleton column selection is random sampling, wherein indices are drawn from a probability distribution that is designed to bias towards structurally important columns. Various sampling distributions have been proposed, including ones that bias toward columns having large norm \cite{frieze_monte_carlo}, or column subsets spanning a large-volume parallelepiped \cite{deshpande_rademacher_volumesampling, deshpande_rademacher_projectiveclustering}. A particularly well-studied sampling strategy, based on the work of Mahoney, Drineas, and Muthukrishnan \cite{mahoney_drineas_muthukrishnan, mahoney_cur}, uses a sampling distribution derived from column-based subspace leverage scores. Given a singular value decomposition $\mA = \mU\mSigma\mV\tp$ and a target rank $k$, the column-based subspace leverage scores are the quantities $\ell_1^{(k)},\, \ldots,\, \ell_n^{(k)} \in [0,\, 1]$ defined by $\ell_j^{(k)} = \| \mV_{j,\, 1:k} \|_2$. These numbers are uniquely defined by $\mA$ (as opposed to depending on the particular $\mV$ used for the SVD) provided that $\sigma_k(\mA) > \sigma_{k + 1}(\mA)$. Going forward, we will simply call these ``leverage scores'' for brevity.

Because $\sum_{j = 1}^n (\ell_j^{(k)})^2 = \| \mV_{:,\, 1:k} \|_\frob^2 = k$, the numbers $P_j = k^{-1} (\ell_j^{(k)})^2$ define a probability distribution over the columns of $\mA$, called the leverage score distribution. Mahoney et al.\ \cite{mahoney_drineas_muthukrishnan, mahoney_cur} have developed an algorithm which select skeleton columns using random draws from $P_j$. Computing leverage scores exactly, which requires computation of singular vectors, can be prohibitively expensive, and therefore, leverage score estimation methods based on random sketching \cite{drineas_fast_approximation} can be instrumental.

As originally described, their algorithm oversamples greatly, returning an approximation whose rank is significantly higher than $k$. We will use a variant of their algorithm which selects $k + p$ column indices by sampling from the leverage score distribution without replacement, where $p$ is a fixed oversampling parameter, and then projects its input column-wise into the subspace spanned by the leading $k$ singular vectors of the skeleton columns. This final projection step is to facilitate a fair comparison with other algorithms that reduce to rank $k$. We will refer to this procedure as leverage score sampling, or (LEVG), and we summarize it in \cref{alg:levg}.
\begin{algorithm}
\caption{\textsc{Leverage Score Sampling (LEVG)}} \label{alg:levg}
\begin{algorithmic}[1]
\STATE $\mU,\, \mSigma,\, \mV \leftarrow \texttt{svd}(\mA)$.
\STATE $P_j \leftarrow k^{-1}\|\mV_{j,\, 1:k} \|_2^2$ for $j = 1,\, \ldots,\, n$.
\STATE $J = (j_1,\, \ldots,\, j_{k + p}) \leftarrow$ a random sample from $P_j$ with replacement.
\STATE $\mQ \leftarrow \texttt{thin\_qr}(\mA_{:,J}),\, \mX \leftarrow \mQ\tp\mA$.
\STATE $\mW, \mGamma,\, \mZ \leftarrow \texttt{svd}(\mX)$.
\STATE $\mV \leftarrow \mQ\mW_{:,\, 1:k}$.
\RETURN $\mB_1 = \mV,\, \mB_2 = \mV\tp\mA$.
\end{algorithmic}
\end{algorithm}

For a more complete overview of interpolative and CUR decompositions, we refer the reader to \cite{murray2023randomized}.

\section{The Randomized Golub-Klema-Stewart Algorithm} \label{section:rgks}

The experiments and analysis in this paper center around a randomizied variant of the GKS algorithm \cite{golub_klema_stewart}, where randomization is incorporated by estimating $\mV_k\tp$ with a randomized singular value decomposition; we call this algorithm RGKS. A related method was introduced by Saibaba \cite{saibaba_rdeim}, who applies a pivoted QR factorization to the left singular vector estimates from the RSVD as part of a model order reduction algorithm called RDEIM. This paper differs from Saibaba's in several ways. First, while RGKS computes $\mA \approx \mP\mA$ where $\mP = \mA_{:,J}(\mA_{:,J})\pinv$ is an orthogonal projection operator, RDEIM uses an oblique projection. Second, RDEIM estimates singular vectors by orthogonalizing the sketch $\mOmega\tp\mA$, where RGKS (via the RSVD) orthogonalizes $\mA\tp\mA\Omega$, thus attaining higher accuracy due to the faster singular value decay in $\mA\tp\mA$ \cite{saibaba_randomized_subspace_iteration}. Finally, Saibaba analyzes RDEIM in the context of approximating a nonlinear, vector-valued function from a small number of its components, and bounds error in terms of additive distance from a ``baseline'' POD mode truncation error. We analyze RGKS in a low-rank approximation context, and bound error in terms of multiplicative distance from the optimal error implied by the Eckart-Young theorem.

Pseudocode for RGKS is provided in \cref{alg:rgks}, where $\texttt{rsvd}(\cdot,\, k,\, p,\, q)$ denotes an RSVD at rank $k$ with oversampling $p$ and $q$ power iterations. For reference, we also provide pseudocode for GKS in \cref{alg:gks}, where $\texttt{partial\_svd}(\cdot,\, k)$ denotes a routine that computes the first $k$ components of an SVD (e.g., by truncating a full SVD or via Krylov methods). Note that when no oversampling is used ($p = 0$), it is possible to replace line \eqref{line:rgks_rsvd} of RGKS with the instructions
\begin{enumerate}
\item $\mOmega \leftarrow$ an $n \times k$ matrix with i.i.d. standard Gaussian entries,
\item $\mQ_1,\, \mR_1 \leftarrow \texttt{thin\_qr}(\mA\mOmega)$,
\item $\mQ_2,\, \mR_2 \leftarrow \texttt{thin\_qr}(\mA\tp\mQ_1)$,
\end{enumerate}
and to execute the rest of the algorithm using $\mQ_2$ in place of $\mV_k$. If one were to compute $\widehat{\mU},\, \widehat{\mSigma}_k,\, \widehat{\mV}_k \leftarrow \texttt{rsvd}(A,\, k,\, 0,\, 0)$, reusing the same $\mOmega$ for the sketch, then $\mQ_2$ and $\widehat{\mV}_k$ would differ only by an orthogonal transformation from the right, meaning the pivoting step in line \eqref{line:rgks_rrqr} of RGKS would select the exact same columns\footnote{We have assumed no power iteration ($q = 0$) to simplify this explanation, but one can extend to $q \geq 1$ by inserting ``$\mQ_1,\, \mR_1 \leftarrow \texttt{thin\_qr}(\mA\mA\tp\mQ_1)$ for $i = 1,\, \ldots,\, q$'' before the final line.}. However, if these instructions are run with $p \geq 1$ (replacing $\mOmega$ with an $n \times (k + p)$ matrix), then it will be necessary to separate the leading $k$-dimensional singular subspace of $\mA\tp\mQ_1$ from the rest of its range space. In this case, computing a rank-revealing factorization of some sort becomes unavoidable, with the SVD being the most natural choice.
\begin{algorithm}
\caption{\textsc{Golub-Klema-Stewart (GKS)}} \label{alg:gks}
\begin{algorithmic}[1]
\STATE $\mU_k,\, \mSigma_k,\, \mV_k \leftarrow \texttt{partial\_svd}(\mA,\, k)$. \label{line:gks_svd}
\STATE $\mPi,\, \mQ,\, \mR \leftarrow \texttt{rrqr}(\mV_k\tp,\, k)$.\label{line:gks_rrqr}
\RETURN $\mB_1 = \mA\mPi_{:,\, 1:k}$, $\mB_2\tp = (\mA\mPi_{:,\, 1:k})\pinv \mA$.
\end{algorithmic}
\end{algorithm}

\begin{algorithm}
\caption{\textsc{Randomized Golub-Klema-Stewart (RGKS)}} \label{alg:rgks}
\begin{algorithmic}[1]
\STATE $\widehat{\mU}_k,\, \widehat{\mSigma}_k,\, \widehat{\mV}_k \leftarrow \texttt{rsvd}(\mA,\, k,\, p,\, q)$. \label{line:rgks_rsvd}
\STATE $\mPi,\, \mQ,\, \mR \leftarrow \texttt{rrqr}(\widehat{\mV}_k\tp,\, k)$. \label{line:rgks_rrqr}
\RETURN $\mB_1 = \mA\mPi_{:,\, 1:k}$, $\mB_2\tp = (\mA\mPi_{:,\, 1:k})\pinv \mA$.
\end{algorithmic}
\end{algorithm}

RGKS is conceptually similar to GRID, in that it selects columns using a pivoted QR factorization of a short matrix obtained through sketching, but the use of sketched singular vectors rather than sketched matrix columns creates a significant difference in behavior, as explained in \cref{section:numerics}. RGKS also bears similarities to leverage score sampling with randomized leverage score estimation, as in \cite{drineas_fast_approximation}, but the use of a pivoted QR rather than random sampling gives RGKS an advantage in terms of accuracy; this is also explained in \cref{section:numerics}. In all of our numerical experiments with GKS and RGKS, we will compute RRQR factorizations using the Golub-Businger algorithm \cite{businger_linear_least_squares}. Many other RRQR factorizations algorithms could be used, e.g., \cite{chandrasekaran_1994, gu_srrqr}.

RGKS has a number of distinctive features that motivate our focus on it. Like GKS, RGKS uses an RRQR factorization to approximately optimize error bounds such as those in \cref{eqn:prior_sigmamin_bounds}, making it highly amenable to error analysis. Second, literature exists which analyzes the accuracy of the singular subspace estimates computed by RSVD, see, e.g., \cite{dong_martinsson_nakatsukasa, saibaba_randomized_subspace_iteration, zhang_tang}. This allows us to explore how the particular subspace errors committed by RSVD affect the accuracy of RGKS. Additionally, in \cref{section:numerics} we will see that the structure of a matrix's right singular subspace (especially its coherence) plays a part in determining the accuracy of many low-rank approximation algorithms. Because RGKS operates directly on an estimate of the right singular subspace, it is a natural case study on the effects of singular subspace structure.   

RGKS is also efficient. Numerical experiments in the supplementary text show that for a fixed value of $k$, RGKS has essentially the same runtime and asymptotic scaling as RSVD, and is only slightly slower than the GRID algorithm \cite{martinsson_rokhlin_tygert_randid, woolfe_srft}.

\section{Numerical Comparison of Algorithms} \label{section:numerics}

\begin{figure}[t!]
    \centering
    \includegraphics[scale=0.7]{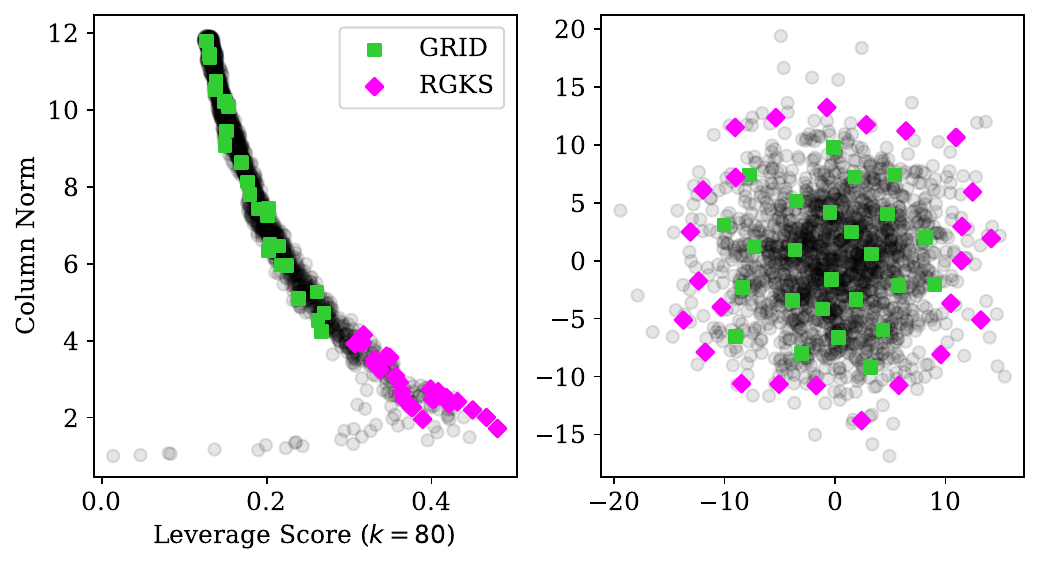}
    \caption{More details on the test matrix $\mA$ from \cref{fig:kahan_kernel_example} (right panel). Left panel: each column of $\mA$ plotted based on its column norm and subspace leverage score. Right panel: the points $x_1,\, \ldots,\, x_{800}$ from which $\mA$ is generated. The first 25 indices selected by GRID and RGKS (from a single trial with $k = 80$) are highlighted.}
    \label{fig:kernel_scatter}
\end{figure}

We begin this section with a set of numerical examples designed to highlight how different interpolative decomposition algorithms respond to different matrix structures. To start, recall \cref{fig:kahan_kernel_example}, which showed RGKS attaining smaller approximation error than GRID on the $800 \times 800$ matrix $\mA_{i,\, j} = \kappa(x_i,\, x_j)$, with $\kappa(x,\, y) = \exp(-\frac{1}{8}\| x - y \|_2^2)$ and $x_1,\, \ldots,\, x_{800}$ drawn from a 2D Gaussian. \Cref{fig:kernel_scatter} shows that the columns of $\mA$ with highest leverage also have some of the smallest norms. Roughly, high-leverage small-norm columns correspond to $x_i$ on the fringes of the distribution, while low-leverage high-norm columns are in the central cluster.

GRID, being a randomized approximation to a true CPQR\footnote{We assume here that the Golub-Businger algorithm is being used inside GRID and RGKS.} of $\mA$'s columns, uses a selection strategy that's biased towards columns of large norm. RGKS, on the other hand, randomly approximates a CPQR of $\mV_k\tp$, whose column norms are precisely the leverage scores, therefore biasing towards columns of high leverage. Bounds in \cref{section:analysis} show that interpolative decomposition error depends on the singular values of $\mV_{J,\, 1:k}$ where $J$ is the skeleton index subset (cf.\ \cref{thm:idbound_subspaceonly,lemma:angle_leverage_connection,thm:idbound_bothstructures}); by approximating a CPQR of $\mV_k\tp$, RGKS is maximizing these singular values more directly than GRID. In \cref{fig:kernel_scatter}, we see how this translates into a qualitative difference between the behaviors of RGKS and GRID. The first 25 columns chosen by GRID are all taken from the central cluster of the Gaussian sample, while the first 25 chosen by RGKS are all taken from the fringes. A qualitative difference is also present in the approximation errors, for as shown in \cref{fig:kahan_kernel_example}, RGKS significantly outperforms GRID.

\begin{figure}
    \includegraphics[scale=.7]{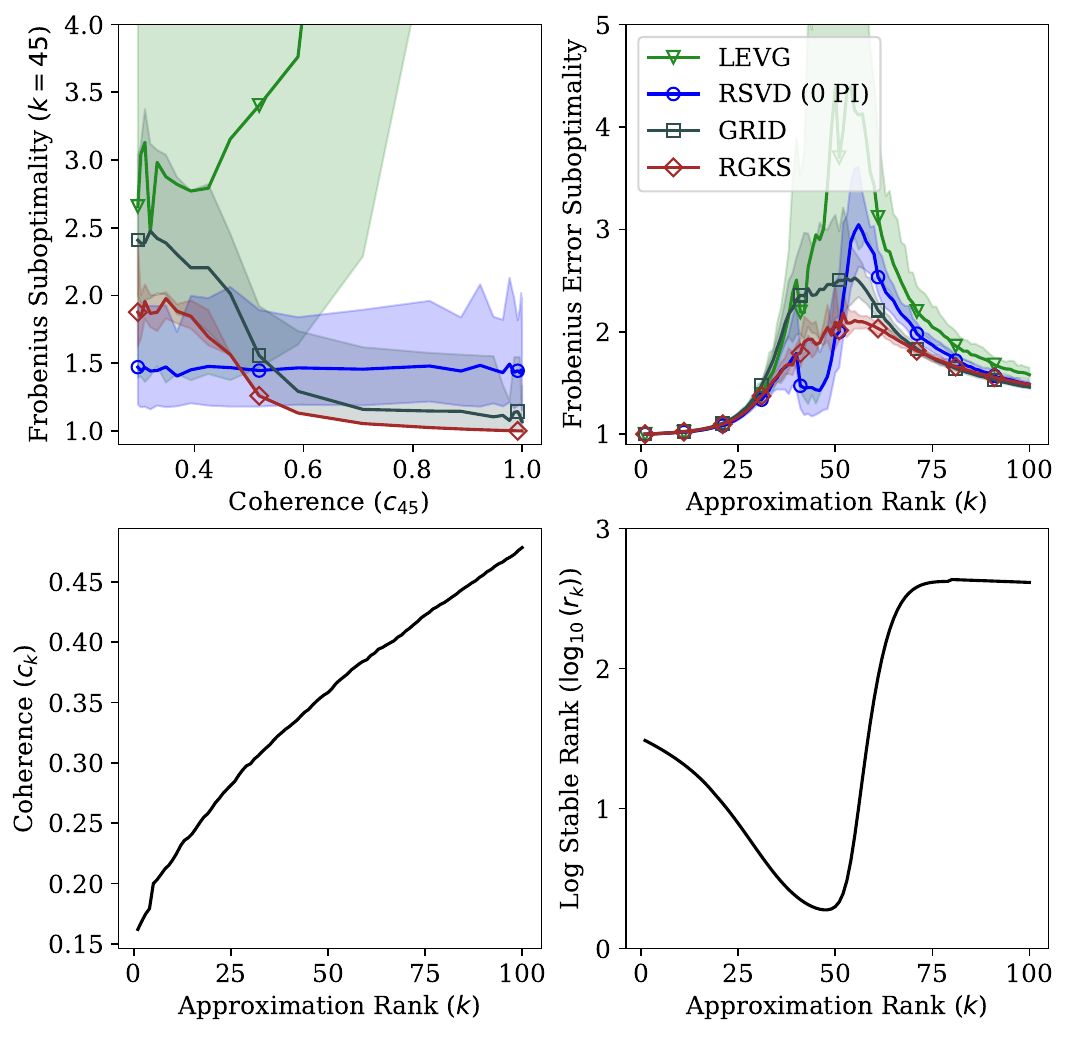}
    \caption{Frobenius norm errors of low-rank approximation algorithms over variations in approximation rank and coherence for $512 \times 512$ test matrices with a fixed singular spectrum. Solid lines indicate a mean over 100 trials, and shaded regions span from the 5\% quantile to the 95\% quantile of the error distribution. The bottom-left panel shows the coherence profile of the single test matrix used for tests over variations in the approximation rank.} \label{fig:crossections}
\end{figure}

\Cref{fig:crossections} illustrates the results of an experiment which, firstly, highlight an important performance difference between RGKS and LEVG. For this experiments, an $n \times n$ test matrix $\mA$ was generated (with $n = 512$) by setting $\mA = \mU\mSigma\mV\tp$, where $\mU$ was an orthogonalization of a Gaussian random matrix. The matrix $\mSigma$ was constructed using a prescribed singular value decay profile, and the coherence of $\mA$'s right singular subspaces was controlled through the construction of $\mV$. Using a permutation matrix for $\mV$ caused $\mA$ to exhibit perfectly coherent right singular subspaces, while using a normalized Hadamard matrix for $\mV$ resulted in perfectly incoherent singular subspaces for $\mA$. To explore intermediate values of coherence,  we computed polar factorizations of convex combinations of a permutation matrix and a normalized Hadamard matrix with varying weights. For each weight, we set $\mV$ to be the orthogonal factor from the polar decomposition, which is the closest orthogonal matrix to the original convex combination \cite[Theorem 8.4]{higham_functions_of_matrices}.

In the bottom row of \cref{fig:crossections}, we have plotted the approximation errors of RSVD, RGKS, GRID, and LEVG over coherence of the leading right singular subspace. While RSVD is completely insensitive to changes in coherence (as predicted by theory), the other interpolative decompositions are much more sensitive to it. A particularly vast performance difference exists between RGKS and LEVG. To understand this, it is informative to examine LEVG without oversampling to facilitate a more direct comparison with the GKS algorithm, which does not oversample during skeleton column selection\footnote{To facilitate a direct comparison, we could also introduce oversampling to GKS by replacing $k$ with $k + p$ in \cref{alg:gks}, followed by an optimal reduction to rank-$k$ in the style of \cref{alg:levg}.}. To that end, consider the matrix
\begin{equation*}
\mM = \begin{bmatrix}
1 & 0 & 1 & 0 \\
0 & 1 & 0 & 1 \\
1 & 0 & 1 & 0 \\
0 & 1 & 0 & 1
\end{bmatrix} = \begin{bmatrix}
\frac{1}{\sqrt{2}} & 0 \\
0 & \frac{1}{\sqrt{2}} \\
\frac{1}{\sqrt{2}} & 0 \\
0 & \frac{1}{\sqrt{2}}
\end{bmatrix} \begin{bmatrix}
2 & 0 \\
0 & 2
\end{bmatrix} \begin{bmatrix}
\frac{1}{\sqrt{2}} & 0 & \frac{1}{\sqrt{2}} & 0 \\
0 & \frac{1}{\sqrt{2}} & 0 & \frac{1}{\sqrt{2}}
\end{bmatrix}.
\end{equation*}
For $k = 2$, the leverage score distribution is uniform over $\mM$'s columns, meaning that LEVG has a 25\% change of picking $J = \{ 1,\, 3 \}$ or $J = \{ 2,\, 4 \}$; for these choices, $\mV_{J,\, 1:k}$ is singular. In contrast, the CPQR factorization in GKS ensures that the selected column indices $J$ satisfy $\sigma_\mathrm{min}(\mV_{J,\, 1:k}) > 0$. While a perfectly uniform leverage score distribution is unlikely to appear in applications, this example nonetheless highlights that LEVG requires a higher degree of oversampling to effectively optimize the bounds in this paper's error analysis; see \cref{thm:idbound_subspaceonly,lemma:angle_leverage_connection,thm:idbound_bothstructures}. To extend these principles to RGKS, one must analyze how $\sigma_\mathrm{min}(\mV_{J,\, 1:k})$ is perturbed under the randomization errors introduced by RSVD; this is the subject of \cref{thm:idbound_randomized_normwise}.

\Cref{fig:crossections} also shows how the performance of RSVD, RGKS, and LEVG (relative to optimality) varies as a function of the approximation rank. The performance comparison between RSVD and RGKS is strongly affected by singular spectrum decay; for example, for the rapidly decaying singular spectrum around $k = 50$, RSVD outperforms RGKS, whereas the behavior is reversed when the decay is slower. This aligns with theoretical analyses \cite{halko_finding_structure_with_randomness, saibaba_randomized_subspace_iteration} which show that RSVD attains a nearly optimal approximation when the singular value gap is large enough. Notice that the RGKS error is, to a very rough approximation, inversely proportional to the residual stable rank $r_k$; \cref{thm:idbound_bothstructures} will formalize this observation in an error bound. An interesting point of comparison is the corresponding spectral norm plot in the supplementary text, which shows that the dependence on $r_k$ is different when errors are measured in the spectral norm.

\section{Structure-Aware Analyses of Interpolative Decompositions} \label{section:analysis} 

Here we develop and analyze error bounds which capture the effects of singular subspace geometry (\cref{subsec:idbound_subspaceonly}) and residual stable rank (\cref{subsec:idbound_bothstructures}), and which exhibit interesting dependence on singular value gaps (\cref{subsec:idbound_spectrumonly}).

\subsection{The Effects of Subspace Geometry} \label{subsec:idbound_subspaceonly}

To develop an interpolative decomposition error bound which captures the effects of subspace geometry, the main idea is to view interpolative decomposition as a process of approximating one subspace by another, namely, approximating $\range(\mA)$ by the span of a small number of $\mA$'s columns. By connecting this to a different problem which involves approximating the leading right singular subspace, we will see that the geometry of this subspace plays an important role in determining the accuracy of the original interpolative decomposition. To work in this setting, we assume that $\sigma_{k}(\mA) > \sigma_{k + 1}(\mA)$ so that the leading right singular subspace $\cV_k$ is well defined. Given this assumption, we define a subspace $\cI_J \subseteq \R^n$ encoding the choice of skeleton columns as
\begin{equation}
\cI_J = \vspan \{ \ve_j \suchthat j \in J \}. \label{eqn:index_support_subspace}
\end{equation}
We now state the main result of this section below as \cref{thm:idbound_subspaceonly}.
\begin{theorem} \label{thm:idbound_subspaceonly}
Choose $k \leq n/2$ such that $\sigma_k(\mA) > \sigma_{k + 1}(\mA)$, let $J = \{ j_1,\, \ldots,\, j_k \}$ be a selection of skeleton column indices, and let $\varphi_\mathrm{max}$ be the largest principal angle between $\cI_J$ and $\cV_k$. If $\varphi_\mathrm{max} < \pi / 2$, then
\begin{equation*}
\| \mA - \mA_{:,\, J}(\mA_{:\, J})\pinv\mA \|_2 \leq \| \mSigma_\perp \|_2 \sec \varphi_\mathrm{max}.
\end{equation*}
\end{theorem}

\begin{proof}
As explained below, this is an equivalent statement of \cite[Theorem 1.5]{hong_pan}. For a self-contained proof, refer to \cref{subsec:idbound_subspace_proof}.
\end{proof}

\Cref{thm:idbound_subspaceonly} shows that the problem of choosing columns of $\mA$ whose span approximates $\range(\mA)$ is related to the problem of choosing elementary unit vectors whose span approximates the right singular subspace $\cV_k$. Leverage scores provide some amount of information on this problem, in that $\ell_j^{(k)}$ measures the degree to which $\ve_j$ is aligned with $\cV_k,$ as seen from the relation
\begin{equation*}
\ell_j^{(k)} = \| \mV_{j,\,1:k} \|_2 = \| \mV_{j,\, 1:k} (\mV_{:,\, 1:k})\tp \|_2 = \| \ve_j\tp\mV_k\mV_k\tp \|_2.
\end{equation*}
\Cref{lemma:angle_leverage_connection} allows us to formalize the connections between leverage scores and $\varphi_\mathrm{max}$.
\begin{lemma} \label{lemma:angle_leverage_connection}
Choose $k \leq n/2$ with $\sigma_k(\mA) > \sigma_{k + 1}(\mA)$, and let $J = \{ j_1,\, \ldots,\, j_k \}$ be a selection of skeleton column indices. If $\varphi_1 \leq \ldots \leq \varphi_k$ are the principal angles between $\cI_J$ and $\cV_k$, then
\begin{equation*}
\cos \varphi_i = \sigma_i(\mV_{J,\, 1:k})
\end{equation*}
for $1 \leq i \leq k$. Furthermore, $\mV_{J,\, 1:k}$ is invertible if and only if $\max_i \varphi_i < \pi / 2$, in which case
\begin{equation*}
\tan \varphi_{k - i + 1} = \sigma_i(\mV_{[n] \setminus J,\, 1:k}(\mV_{J,\, 1:k})^{-1})
\end{equation*}
for $1 \leq i \leq k$.
\end{lemma}
\begin{proof}
Refer to \cref{subsec:angle_svd_connection_proof}.
\end{proof}
An important consequence of \cref{lemma:angle_leverage_connection} is that $\cos \varphi_\mathrm{max} = \sigma_\mathrm{min}(\mV_{J,\, 1:k})$. This, together with \cref{eqn:general_cpqr_id}, proves the equivalence of \cref{thm:idbound_subspaceonly} and \cite[Theorem 1.5]{hong_pan}.

\Cref{lemma:angle_leverage_connection} provides a means of quantifying the effects of subspace geometry on interpolative decomposition accuracy, because it implies upper and lower bounds on $\sec \varphi_\mathrm{max}$ that depend on leverage scores and coherence. Recall that coherence at rank $k$ is defined as $c_k = \max_{1 \leq j \leq n} \ell_j^{(k)}$. In the case of low coherence, we can place a lower bound on $\sec \varphi_\mathrm{max}$ by using \cref{lemma:angle_leverage_connection}, together with the fact that the minimum singular value of a matrix cannot exceed the minimum row norm:
\begin{equation*}
\sec \varphi_\mathrm{max} = \frac{1}{\sigma_\mathrm{min}(\mV_{J,\, 1:k})} \geq\frac{1}{\min_{j \in J} \ell_j^{(k)}} \geq c_k^{-1}.
\end{equation*}
In the case of near-minimal coherence, we have $c_k \leq (1 + \varepsilon)\sqrt{n^{-1}k}$ for some small $\varepsilon > 0$. Then $\sec \varphi_\mathrm{max} \geq (1 + \varepsilon)^{-1}\sqrt{k^{-1}n}$, and as expected for the incoherent case, this bound holds independent of the choice of $J$. 

For the case where there are $k$ leverage scores with values near unity and corresponding to column indices $J = \{ j_1,\, \ldots,\, j_k \}$, we can bound $\sec \varphi_\mathrm{max}$ from above. Using \cref{lemma:angle_leverage_connection},
\begin{equation*}
\sum_{i = 1}^k \cos^2 \varphi_i = \sum_{i = 1}^k \sigma_i(\mV_{J,\, 1:k})^2 = \| \mV_{J,\, 1:k} \|_\frob^2 = \sum_{j \in J} (\ell_j^{(k)})^2.
\end{equation*}
This implies $\cos^2 \varphi_\mathrm{max} = {\textstyle \sum_{j \in J} (\ell_j^{(k)})^2 - \sum_{i \neq i_\mathrm{max}} \cos^2 \varphi_i}$, where $i_\mathrm{max} \in \argmax_i \varphi_i$. Inserting $\cos \varphi_i \leq 1$ for $i \neq i_\mathrm{max}$, and assuming the skeleton columns have large enough leverage scores that $\sum_{j \in J} (\ell_j^{(k)})^2 \geq k - 1$, we obtain
\begin{equation}
\sec \varphi_\mathrm{max} \leq \frac{1}{\sqrt{1 - \sum_{j \in J} (1 - (\ell_j^{(k)})^2)}}, \label{eqn:highcoherence_secantbound}
\end{equation}
a bound which is small when $\ell_j^{(k)} \approx 1$ for all $j \in J$. Choosing $k$ indices with leverage scores near unity is sufficient, but emphatically \emph{not} necessary, for $\sec \varphi_\mathrm{max}$ to be small. For example, \cref{fig:crossections} shows that even at the minimum of coherence, when there are no leverage scores near unity, RGKS can still attain close to optimal error.

Although bounds on $\sec \varphi_\mathrm{max}$ are useful in describing the effects of coherence and leverage scores, $\sec \varphi_\mathrm{max}$ itself can be overly pessimistic as a bound on approximation error (relative to $\| \mSigma_\perp \|_2$). To see this, consider any subset $I \subseteq J$ of size $k - t$, and assuming that $\sigma_{k - t}(\mA) > \sigma_{k - t + 1}(\mA)$, let $\varphi_\mathrm{max}\ps{I}$ denote the largest principal angle between $\cV_{k - t}$ and $\cI_I$. Using \cref{thm:idbound_subspaceonly}, together with the fact that interpolative decomposition error decreases as skeleton columns are added,
\begin{align}
\| \mA - \mA_{:,\, J}(\mA_{:,\, J})\pinv\mA \|_2 &\leq \| \mA - \mA_{:,\, I}(\mA_{:,\, I})\pinv\mA \|_2 \nonumber \\
&\leq \sigma_{k - t + 1}(\mA)\sec \varphi_\mathrm{max}\ps{I} \nonumber \\
&= \frac{\| \mSigma_\perp \|_2 \sec \varphi_\mathrm{max}\ps{I}}{\gamma_{k - t + 1,\, k + 1}}, \label{eqn:idbound_subspaceonly_subsetangle}
\end{align}
where $\gamma_{i,\, j} = \sigma_{j}(\mA)/\sigma_{i}(\mA)$. Equivalently, because of \cref{lemma:angle_leverage_connection},
\begin{equation}
\| \mA - \mA_{:,\, J}(\mA_{:,\, J})\pinv\mA \|_2 \leq \frac{\| \mSigma_\perp \|_2}{\gamma_{k - t + 1,\, k + 1} \sigma_\mathrm{min}(\mV_{I,\, 1:k - t})}. \label{eqn:idbound_subspaceonly_subsetsigma}
\end{equation}

Equation \cref{eqn:idbound_subspaceonly_subsetangle} shows that when $\sigma_{k - t + 1}(\mA) \approx \sigma_{k + 1}(\mA)$, obtaining a near-optimal approximation does not necessarily require that $\varphi_\mathrm{max}$ is small. Rather, it is sufficient that a smaller singular subspace $\cV_{k - t}$ is well-approximated by some smaller subset of elementary unit vectors, corresponding to a small value of $\varphi_\mathrm{max}\ps{I}$. From the perspective of minimal singular values, equation \cref{eqn:idbound_subspaceonly_subsetsigma} shows that a large value of $\sigma_\mathrm{min}(\mV_{J,\, 1:k})$ is not strictly necessary for obtaining near-optimal approximation error. Instead, it is sufficient that $\mV_{J,\, 1:k}$ contains a submatrix $\mV_{I,\, 1:k - t}$ whose minimal singular value is large, for some $t$ such that $\sigma_{k - t + 1}(\mA) \approx \sigma_{k + 1}(\mA)$. This is easiest to achieve when $t$ is chosen as large as possible while maintaining $\sigma_{k - t + 1}(\mA) \approx \sigma_{k + 1}(\mA)$.

\begin{figure}
    \includegraphics[scale=0.63]{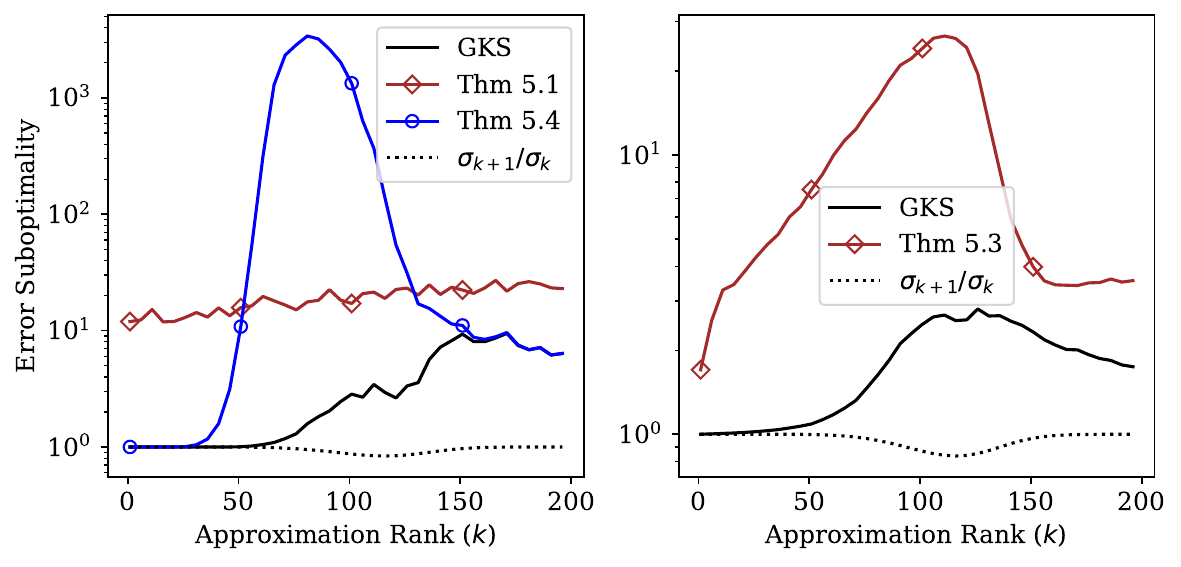}
    \caption{Errors of GKS relative to optimality in the spectral (left) and Frobenius (right) norms compared with error bounds, versus approximation rank, for the same $512 \times 512$ test matrix as in \cref{fig:crossections}.} \label{fig:errorbounds_rank}
\end{figure}

The left panel of \cref{fig:errorbounds_rank} compares the error bound in \cref{thm:idbound_subspaceonly} to the actual error for approximations computed with RGKS. In the right-hand plot, we see that the error bound accurately reflects the changes in approximation accuracy as subspace geometry (measured by coherence) is varied. However, note how the looseness of $\sec \varphi_\mathrm{max}$ as a bound on approximation error (relative to $\| \mSigma_\perp \|_2$) is illustrated in the left-hand plot, particularly for small approximation ranks.
\begin{figure}
    \includegraphics[scale=0.63]{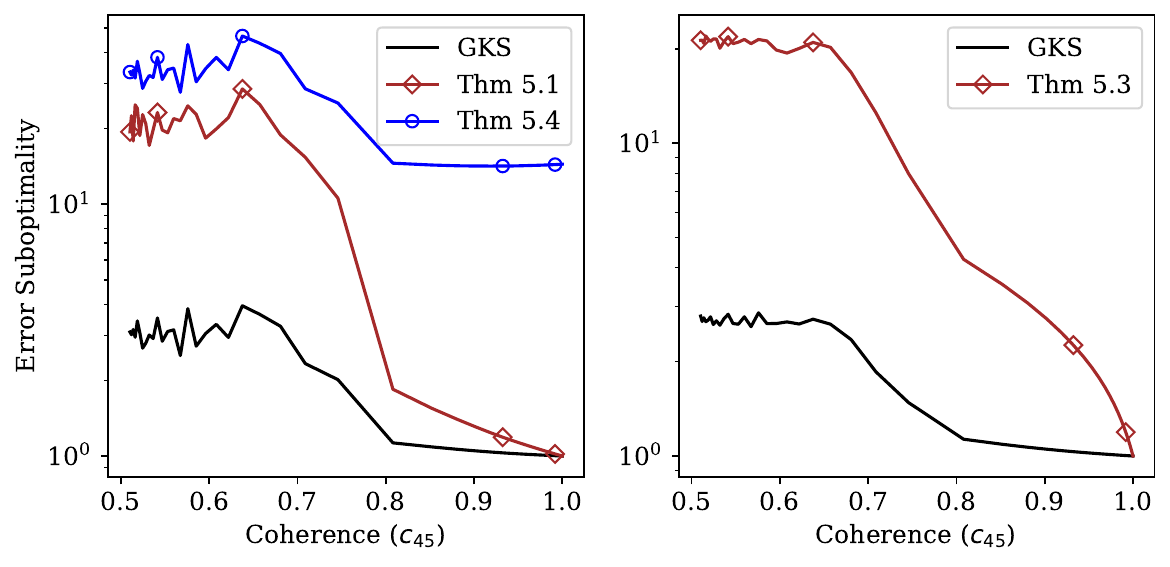}
    \caption{Errors of GKS relative to optimality in the spectral (left) and Frobenius (right) norms at $k = 45$ compared with error bounds, versus coherence. All test matrices were $512 \times 512$, and had the same singular spectrum as the matrices in \cref{fig:crossections,fig:errorbounds_rank}.} \label{fig:errorbounds_coher}
\end{figure}

\subsection{Subspace Geometry and Stable Rank} \label{subsec:idbound_bothstructures}

The bound in the previous section captures the effects of subspace structure, but fails to capture the effects of spectral decay structures such as singular value gap and stable rank. The left panel of \cref{fig:errorbounds_rank}, for example, shows that the relative spectral error of RGKS tends to be greatest in or near regions with both high stable rank and a significant singular value gap, while the plot of the corresponding error bound shows no such dependence on the singular spectrum. We now present a bound which captures the combined effects of subspace geometry and singular spectral decay, where spectral decay is measured by the residual stable rank; this is accomplished by shifting our analysis to errors measured in the Frobenius norm. We state the result below in \cref{thm:idbound_bothstructures}.
\begin{theorem} \label{thm:idbound_bothstructures}
Choose $k \leq n/2$ such that $\sigma_k(\mA) > \sigma_{k + 1}(\mA)$, let $J = \{ j_1,\, \ldots,\, j_k \}$ be a selection of skeleton column indices, and let $\varphi_1,\, \ldots,\, \varphi_k$ be the principal angles between $\cI_J$ and $\cV_k$. If $\max_i \varphi_i < \pi / 2$, then
\begin{equation}
\| \mA - \mA_{:,\, J}(\mA_{:,\, J})\pinv\mA \|_\frob \leq \| \mSigma_\perp \|_\frob \sqrt{ 1 + \frac{1}{r_k} \sum_{i = 1}^k \tan^2 \varphi_i }, \label{eqn:idbound_bothstructures_allangles}
\end{equation}
where $r_k = \| \mSigma_\perp \|_\frob^2 / \| \mSigma_\perp \|_2^2$ is the residual stable rank.
\end{theorem}

\begin{proof}
Refer to \cref{subsec:idbound_bothstructures_proof}.
\end{proof}

\Cref{thm:idbound_bothstructures} depends on subspace geometry effects through $\tan \varphi_1,\, \ldots,\, \tan \varphi_k$, and depends on all of the principal angles between $\cI_J$ and $\cV_k$. (Whereas the bound from \cref{subsec:idbound_subspaceonly} depended on only the largest angle.) Unlike the previous bound, \cref{thm:idbound_bothstructures} also incorporates the effects of stable rank in the form of a factor $r_k^{-1}$ that attenuates the impact of subspace geometry. This suggests that when errors are measured in the Frobenius norm, a near-optimal interpolative decomposition is easier to obtain at approximation ranks corresponding to a high value of $r_k$, an idea which is confirmed by the experiments plotted in \cref{fig:crossections}. Further confirmation is offered by the right panel of \cref{fig:errorbounds_rank}, which plots the actual RGKS error compared with the error bound in \cref{thm:idbound_bothstructures} as the approximation rank and coherence are varied.

\subsection{An Error Bound for Flat Singular Spectra} \label{subsec:idbound_spectrumonly}

In this section we bound the approximation error of interpolative decompositions at ranks where the absence of singular value decay renders \cref{thm:idbound_subspaceonly} and \cref{thm:idbound_bothstructures} inapplicable. In addition, the bound in this section reveals that the relative error of a low-rank approximation is strongly related to the spectral properties of the actual residual error matrix. While it is common for error bounds to depend on the spectrum of the optimal residual singular value matrix $\mSigma_\perp$, the bound in this section depends on the spectrum of the actual residual matrix.

Consider a subspace $\cW \subseteq \R^m$ with $\dim \cW = k < m$, and define $\mE = \mA - \mP_\cW \mA$, where $\mP_\cW$ is the orthogonal projection matrix onto $\cW$. We can interpret $\mE$ as the residual of a low-rank approximation computed by projecting the columns of $\mA$ into $\cW$. Furthermore, we define the truncated condition number at level $i$ of a matrix $\mM$ as $\kappa(\mM,\, i) \defeq \sigma_1(\mM)/\sigma_i(\mM)$. We can now state \cref{thm:lowrankbound_spectrumonly} below, the main result of this section.
\begin{theorem} \label{thm:lowrankbound_spectrumonly}
Let $\cW \subseteq \R^m$ be a $k$-dimensional subspace, generating a low-rank approximation $\mA \approx \mP_\cW\mA$, and let $\mE = (\mI - \mP_\cW)\mA$ be the error of this approximation. If $\sigma_{k + 1}(\mE) > 0$, then $\| \mE \|_2 \leq \| \mSigma_\perp \|_2 \kappa(\mE,\, k + 1)$.
\end{theorem}

\begin{proof}
Refer to \cref{subsec:lowrankbound_spectrumonly_proof}.
\end{proof}

The left panel of \Cref{fig:errorbounds_rank} illustrates the approximation error of RGKS versus the bound in \cref{thm:lowrankbound_spectrumonly}. At approximation ranks where there is a large singular value gap, the bound is extremely loose. But in certain regions where the singular spectrum does not decay (or decays slowly), the bound is strikingly tight. A result which explains this behavior is provided in the supplementary text.

\section{Analysis of RGKS} \label{section:rgks_analysis}

The error bounds discussed so far are applicable to any interpolative decomposition, regardless of the algorithm used to compute it. If, however, the decomposition is computed by RGKS, then these bounds become especially useful in in that they depend on quantities that are optimized more directly in RGKS than in other interpolative decomposition algorithms. This allows for the development of more refined error bounds for RGKS. \Cref{subsec:gks_bounds} will begin by developing bounds for GKS, where the lack of randomization makes things more straightforward. \Cref{subsec:rgks_extension} then extends the GKS error bounds to RGKS by building off of a preexisting error analysis for RSVD. In practice, we find that RGKS is more robust to randomization errors than the arguments in \cref{subsec:rgks_extension} suggest, and \cref{subsec:rgks_robustness} develops theory which offers partial explanations for this behavior. \Cref{subsec:subspace_approximation_numerics} presents numerical experiments supporting the arguments in \cref{subsec:rgks_robustness}.

\subsection{Analysis of GKS} \label{subsec:gks_bounds}

\Cref{thm:idbound_subspaceonly} bounds the spectral norm error of a general interpolative decomposition in terms of $\sec \varphi_\mathrm{max}$, where $\varphi_\mathrm{max}$ is the largest principal angle between $\cV_k$ and $\cI_J \defeq \vspan\{ \ve_j \suchthat j \in J \}$. Unlike in other interpolative decomposition algorithms, the column selection strategy of GKS can be viewed directly in terms of minimizing $\varphi_\mathrm{max}$. Indeed, if we write the RRQR factorization in line \ref{line:gks_rrqr} of GKS (\cref{alg:gks}) as
\begin{equation*}
\mV_k\tp\mPi = \begin{bmatrix} (\mV_{J,\, 1:k})\tp & (\mV_{[n] \setminus J,\, 1:k})\tp \end{bmatrix} = \mQ \begin{bmatrix} \mR_1 & \mR_2 \end{bmatrix},
\end{equation*}
with $\mR_1 \in \R^{k \times k}$, then \cref{lemma:angle_leverage_connection} implies that $\sec \varphi_\mathrm{max} =  \sigma_\mathrm{min}(\mV_{J,\, 1:k})^{-1} = \sigma_\mathrm{min}(\mR_1)^{-1}$. Rank-revealing QR factorizations are designed explicitly to choose a well conditioned set of basis columns, i.e., a $J$ for which $\sigma_\mathrm{min}(\mR_1)$ is large. In this sense, line \ref{line:gks_rrqr} of GKS serves directly to minimize $\varphi_\mathrm{max}$. Quantitatively, the algorithmic guarantees of RRQR factorizations are such that $\sigma_\mathrm{min}(\mR_1) \geq \sigma_k(\mV_k\tp)q(n,\, k)^{-1} = q(n,\, k)^{-1}$, where $q$ is a function whose growth in $n$ is bounded by a low-degree polynomial. Therefore, \cref{lemma:angle_leverage_connection} shows that $\varphi_\mathrm{max} \leq \arccos(q(n,\, k)^{-1})$, meaning by \cref{thm:idbound_subspaceonly} that
\begin{align}
\| \mA - \mA_{:,\, J}(\mA_{:,\, J})\pinv\mA \|_2 &\leq \| \mSigma_\perp \|_2 \sec(\arccos(q(n,\, k)^{-1})) \label{eqn:gksbound_trigonometric} \\
&\leq  q(n,\, k) \| \mSigma_\perp \|_2. \nonumber
\end{align}
Interestingly, the same error bound holds for interpolative decompositions obtained by applying a RRQR factorization directly to the columns of $\mA$, rather than to the rows of $\mV_k$ as in GKS. This follows from equations \cref{eqn:general_cpqr_id} and \cref{eqn:rrqr_bounds}. The GKS bound can be refined by considering specific RRQR factorization algorithms for which an explicit form of $q(n,\, k)$ is available. For example, if the Gu-Eisenstat algorithm \cite{gu_srrqr} is used in line \ref{line:gks_rrqr} of GKS, then the bound becomes
\begin{equation}
\| \mA - \mA_{:,\, J}(\mA_{:,\, J})\pinv\mA \|_2 \leq \| \mSigma_\perp \|_2 \sqrt{1 + f^2k(n - k)}, \label{eqn:gks_gu_eisenstat_spectral}
\end{equation}
where $f > 1$ is a user-selected parameter controlling a trade-off between maximizing $\sigma_\mathrm{min}(\mR_1)$ and minimizing runtime.

For the Frobenius norm error bound, \cref{thm:idbound_bothstructures}, similar reasoning applies. In addition to terms that depend only on the singular values of $\mA$, this bound depends on skeleton column choices through the term $\sum_i \tan^2 \varphi_i$. \Cref{lemma:angle_leverage_connection} allows us to express this in terms of a quantity that GKS approximately optimizes:
\begin{align*}
\sum_{i = 1}^k \tan^2 \varphi_i &= \sum_{i = 1}^k \sigma_i(\mV_{[n] \setminus J,\, 1:k}(\mV_{J,\, 1:k})^{-1})^2 = \sum_{i = 1}^k \sigma_i(\mR_2\tp \mR_1\ntp)^2 \\
&= \| \mR_1^{-1}\mR_2 \|_\frob^2.
\end{align*}
The quantity $\| \mR_1^{-1}\mR_2 \|_\frob^2$ measures the size of the coefficients used by the RRQR factorization to interpolate the rows of $\mV_k$ in terms of $\mV_{J,\, 1:k}$. The RRQR factorization in GKS strives to minimize these coefficient magnitudes by choosing a well-conditioned row basis, and in some cases, upper-bounds are available. For example, the Gu-Eisenstat algorithm chooses $J$ such that the entries of $\mR_1^{-1}\mR_2$ are bounded in magnitude by $f$ \cite{gu_srrqr}. This leads to the bound $\| \mR_1^{-1}\mR_2 \|_\frob \leq f \sqrt{k(n - k)}$ and, by an application of \cref{thm:idbound_bothstructures},
\begin{equation}
\| \mA - \mA_{:,\, J}(\mA_{:,\, J})\pinv\mA \|_\frob \leq \| \mSigma_\perp \|_\frob \sqrt{1 + r_k^{-1}f^2 k(n - k)}. \label{eqn:gks_gu_eisenstat_frob}
\end{equation}
Notice the improvement over \cref{eqn:gks_gu_eisenstat_spectral} by a factor of $r_k^{-1}$ inside the square-root. This improvement suggests, in a limited sense, that GKS with the Gu-Eisenstat algorithm is a near-optimal interpolative decomposition for certain classes of matrices. Specifically, it is known \cite[Theorem 3]{deshpande_rademacher_volumesampling} that for any $k \geq 1$ and $\varepsilon > 0$, there exists a $k \times (k + 1)$ matrix $\mA$ such that
\begin{equation}
\min_{K \subseteq [n],\, |K| = k} \| \mA - \mA_{:,\, K}(\mA_{:,\, K})\pinv\mA \|_2 \geq \| \mSigma_\perp \|_\frob (1 - \varepsilon)\sqrt{1 + k}. \label{eqn:deshpande_rademacher_lowerbound}
\end{equation}
Therefore, the smallest interpolative decomposition error bound which can hold over all matrices and all $k$ is $\| \mSigma_\perp \|_\frob \sqrt{1 + k}$. GKS with the Gu-Eisenstat algorithm nearly attains this bound for matrices with large residual stable rank ($r_k \to n - k$).

\subsection{Extension to RGKS} \label{subsec:rgks_extension}

We now extend the analysis to RGKS, where the main challenge is to account for errors introduced by randomizing the SVD. The analysis in this section builds on well-developed theory for RSVD and uses a measure of randomization error that has been extensively studied in previous literature. However, the resulting bounds on RGKS error will only be tight when $k$ corresponds to a large singular value gap. \Cref{subsec:rgks_robustness} will offer two alternative analyses that can better handle shallow singular value gaps, but the first analysis uses a comparison with optimal error at a rank smaller than $k$, while the second analysis uses a measure of randomization error that is less well studied.

The column selection strategy of RGKS (\cref{alg:rgks}) is equivalent to running GKS on the rank-$k$ approximation $\widehat{\mA} = \widehat{\mU}_k\widehat{\mSigma}_k\widehat{\mV}_k\tp$ computed by RSVD in line \ref{line:rgks_rsvd}. Therefore, while GKS approximately minimizes the largest angle between $\cI_J$ and $\cV_k$, RGKS instead minimizes the largest angle between $\cI_J$ and $\widehat{\cV}_k$, where $\widehat{\cV}_k = \range(\widehat{\mV}_k)$ is the RSVD estimate of $\cV_k$. If we denote this largest angle by $\widehat{\varphi}_\mathrm{max}$, then extending the analysis of the previous section means bounding $\varphi_\mathrm{max}$ from above in terms of $\widehat{\varphi}_\mathrm{max}$. A bound of this sort must depend on the accuracy of the singular vector estimates produced by RSVD; this can be accounted for using $\theta_\mathrm{max}$, the largest principal angle between $\cV_k$ and $\widehat{\cV}_k$. Principal angles have long been used as a measure of subspace approximation errors, most famously in the perturbation theory developed by Davis and Kahan \cite{davis_kahan}, Stewart \cite{stewart_1973}, and Wedin \cite{wedin_sine_theta}. Recent years have also seen the developments of principal angle bounds that take into account the algorithmic details of RSVD \cite{dong_martinsson_nakatsukasa, saibaba_randomized_subspace_iteration}.

Using $\theta_\mathrm{max}$ to quantify randomization error, \cref{thm:idbound_randomized_normwise} provides a perturbation bound for extending the GKS analysis to RGKS.
\begin{theorem} \label{thm:idbound_randomized_normwise}
Choose $k \leq n/2$ is such that $\sigma_k(\mA) > \sigma_{k + 1}(\mA)$. Given skeleton indices $J = \{j_1,\, \ldots,\, j_k \}$ and a $k$-dimensional subspace $\widehat{\cV}_k \approx \R^n$, let $\theta_\mathrm{max}$ be the largest principal angle between $\widehat{\cV}_k$ and $\cV_k$, and let $\varphi_\mathrm{max}$ (resp. $\widehat{\varphi}_\mathrm{max}$) be the largest principal angle between $\cI_J$ and $\cV_k$ (resp. $\widehat{\cV}_k$). Then,
\begin{equation*}
\varphi_\mathrm{max} \leq \widehat{\varphi}_\mathrm{max} + \theta_\mathrm{max}.
\end{equation*}
\end{theorem}

\begin{proof}
An elementary proof is given in \cref{subsec:idbound_randomized_proof}. This result also follows from the more general triangle inequalities proven in \cite{qiu_grassmann_2005} for symmetric gauge functions over principal angles.
\end{proof}

Combining \cref{thm:idbound_subspaceonly} with \cref{thm:idbound_randomized_normwise}, we obtain the RGKS error bound
\begin{equation}
\| \mA - \mA_{:,\, J}(\mA_{:,\, J})\pinv\mA \|_2 \leq \| \mSigma_\perp \|_2 \sec (\widehat{\varphi}_\mathrm{max} + \theta_\mathrm{max}). \label{eqn:rgks_spectralbound}
\end{equation}

Line \ref{line:rgks_rrqr} of RGKS ensures that $\widehat{\varphi}_\mathrm{max} \leq \arccos(q(n,\, k)^{-1})$, where the form of $q$ depends on the specific algorithm used to compute the RRQR factorization. In practice, this bound on $\widehat{\varphi}_\mathrm{max}$ may be quite loose. For a Frobenius norm bound, we replace the upper bound in \cref{thm:idbound_bothstructures} with the larger bound $\| \mSigma_\perp \|_2 (1 + kr_k^{-1} \tan^2 \varphi_\mathrm{max})^{1/2}$. Then, applying \cref{thm:idbound_randomized_normwise}, we have
\begin{equation}
\| \mA - \mA_{:,\, J}(\mA_{:,\, J})\pinv\mA \|_\frob \leq \| \mSigma_\perp \|_\frob \sqrt{1 + \frac{k}{r_k}\tan^2 (\widehat{\varphi}_\mathrm{max} + \theta_\mathrm{max})}. \label{eqn:rgks_frobbound}
\end{equation}

\begin{figure}
    \centering
    \includegraphics[scale=.6]{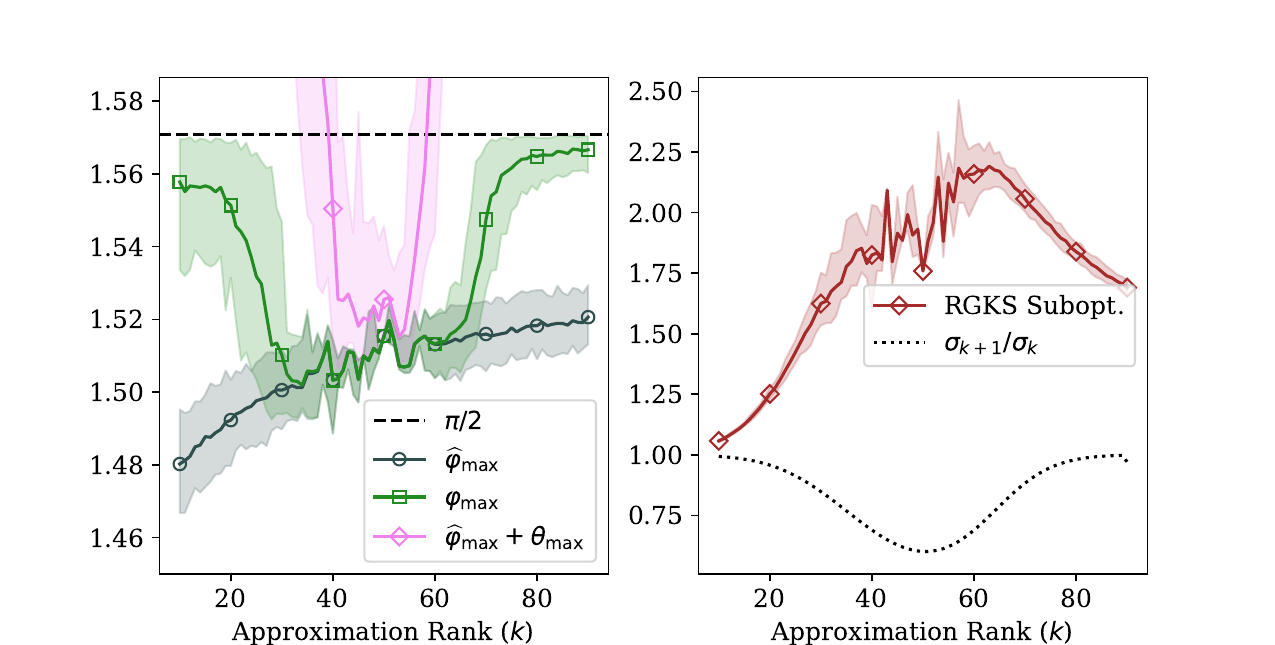}
    \caption{Principal angles and Frobenius norm approximation errors (relative to optimality) for approximations computed by RGKS at different ranks for a $512 \times 512$ test matrix. Each data point is the mean over 100 independent runs of RGKS, with shaded regions indicating 10\% and 90\% quantiles.}
    \label{fig:angle_perturbations}
\end{figure}

The usefulness of these bounds depends on the size of $\theta_\mathrm{max}$. In particular, for these bounds to be non-vacuous, the RSVD in line \ref{line:rgks_rsvd} of RGKS must be accurate enough that $\theta_\mathrm{max} < \pi / 2 - \widehat{\varphi}_\mathrm{max}$. Saibaba \cite{saibaba_randomized_subspace_iteration} has shown that
\begin{equation}
\expect[\sin \theta_\mathrm{max}] \leq \frac{\gamma_k^{2q + 2} C(p)}{\sqrt{1 + \gamma_k^{4q + 4} C(p)^2}}, \label{eqn:saibaba_expectationbounds}
\end{equation}
where $q$ is the RSVD power iteration number, $p$ is the oversampling parameter, and $C$ is a decreasing function of $p$. This result shows that $\theta_\mathrm{max}$ is small with high probability when $\gamma_k \ll 1$, but this is unfortunately not the case when $\gamma_k \approx 1$. The situation is illustrated in \cref{fig:angle_perturbations}, which shows $\varphi_\mathrm{max},\, \widehat{\varphi}_\mathrm{max}$, and $\widehat{\varphi}_\mathrm{max} + \theta_\mathrm{max}$ as a function of $k$ for a test matrix having widely varying values of $\gamma_k$. At $k$ where $\gamma_k$ is very small, $\widehat{\varphi}_\mathrm{max} + \theta_\mathrm{max}$ bounds $\varphi_\mathrm{max}$ well below $\pi / 2$. But as $\gamma_k$ approaches unity, $\widehat{\varphi}_\mathrm{max} + \theta_\mathrm{max}$ rapidly exceeds $\pi / 2$.

\subsection{Robustness of RGKS to Noise} \label{subsec:rgks_robustness}

When $\gamma_k \approx 1$, the actual error of RGKS may be smaller than the previous section's analysis would suggest. For example, consider \cref{fig:angle_perturbations}, and in particular, examine the values of $k$ where $\widehat{\varphi}_\mathrm{max} + \theta_\mathrm{max}$ approaches $\pi / 2$. In this limit, the error bounds in \cref{eqn:rgks_spectralbound} and \cref{eqn:rgks_frobbound} diverge to $+\infty$. However, $\varphi_\mathrm{max}$ itself does not approach $\pi / 2$ at these ranks, and more importantly, the actual RGKS error remains reasonably close to optimal.

One potential explanation for this behavior is that even when the RSVD in line \ref{line:rgks_rsvd} of RGKS estimates $\cV_k$ poorly, a smaller singular subspace $\cV_{k - t}$ might be estimated quite well. To formalize this, assume that $\sigma_{k - t}(\mA) > \sigma_{k - t + 1}(\mA)$ so that $\cV_{k - t}$ is well-defined. Let $\widehat{\cV}_{k - t}$ denote the subspace spanned by the first $k - t$ right singular vectors estimated by the RSVD, and define $\theta_\mathrm{max}\ps{k - t}$ to be the largest principal angle between $\cV_{k - t}$ and $\widehat{\cV}_{k - t}$. Because the RSVD is computed at rank $k$ with oversampling $p$, this construction of $\widehat{\cV}_{k - t}$ is equivalent to running RSVD at rank $k - t$ with oversampling $p + t$. Therefore, \cref{eqn:saibaba_expectationbounds} shows that
\begin{equation}
\expect[\sin \theta_\mathrm{max}\ps{k - t}] \leq \frac{\gamma_{k - t}^{2q + 2} C(p + t)}{\sqrt{1 + \gamma_{k - t}^{4q + 4} C(p + t)^2}}. \label{eqn:saibaba_smaller_subspace}
\end{equation}
If $\gamma_{k - t} \ll 1$, then this bound implies that $\theta_\mathrm{max}\ps{k - t}$ is small with high probability, regardless of the value of $\gamma_k$. This is in contrast to $\theta_\mathrm{max}$, which is essentially uncontrolled when $\gamma_k \approx 1$. The extra oversampling in \cref{eqn:saibaba_smaller_subspace} as compared to \cref{eqn:saibaba_expectationbounds} makes the difference between $\theta_\mathrm{max}\ps{k - t}$ and $\theta_\mathrm{max}$ even more pronounced. A proof in \cref{subsec:rgks_smallerskeleton_proof} shows that for any $I \subseteq J$, RGKS satisfies the error bound
\begin{equation}
\| \mA - \mA_{:,\, J}(\mA_{:,\, J})\pinv\mA \|_2 \leq \frac{\| \mSigma_\perp \|_2}{\gamma_{k - t + 1,\, k + 1}} \sec( \widehat{\varphi}_\mathrm{max}^{(I)} + \theta_\mathrm{max}\ps{k - t}),\label{eqn:rgksbound_smallerskeleton}
\end{equation}
where $\gamma_{k - t + 1,\,k + 1} = \sigma_{k + 1}(\mA)/\sigma_{k - t + 1}(\mA)$, and where $\widehat{\varphi}_\mathrm{max}^{(I)}$ is the largest principal angle between $\cI_I$ and $\cV_{k-t}$. Because of the difference between $\theta_\mathrm{max}$ and $\theta_\mathrm{max}\ps{k - t}$, this bound may not suffer from the divergence to $+\infty$ when the spectrum is flat after $\sigma_k(\mA)$. However, this analysis pays the price of comparing to optimal error at a rank smaller than $k$, as reflected by the presence of $\gamma_{k - t + 1,\,k + 1}$ in \cref{eqn:rgksbound_smallerskeleton}.

Another potential explanation for why RGKS error can remain small, even in the absence of rapid singular value decay, is that the upper bound on $\varphi_\mathrm{max}$ in \cref{thm:idbound_randomized_normwise} may be overly pessimistic given the kinds of subspace errors that RSVD commits. This is reasonable to expect given that $\theta_\mathrm{max}$ is a aggregated measure of error, in the sense that $\sin \theta_\mathrm{max} =  \| \mP - \widehat{\mP} \|_2$, where $\mP$ and $\widehat{\mP}$ are orthogonal projection matrices onto $\cV_k$ and $\widehat{\cV}_k$. The norm of $\mP - \widehat{\mP}$ captures the aggregated total of RSVD subspace approximation error, but it does not capture how the RSVD errors are distributed across the components of the subspace. It may be that component-wise errors, when distributed correctly, perturb $\varphi_\mathrm{max}$ less than \cref{thm:idbound_randomized_normwise} would suggest. Therefore, we will now analyze how $\varphi_\mathrm{max}$ is perturbed under subspace approximation errors with controlled component-wise distributions.

An important measure of component-wise error is the row-wise subspace distance, defined as follows: choose orthonormal bases $\mV_k$ and $\widehat{\mV}_k$ for $\cV_k$ and $\widehat{\cV}_k$, respectively, and for $1 \leq j \leq n$, let $\vv_j$ and $\widehat{\vv}_j$ denote the $j\nth$ rows of $\mV_k$ and $\widehat{\mV}_k$. The row-wise distance between $\cV_k$ and $\widehat{\cV}_k$ is given by\begin{equation}
d_\mathrm{row}(\cV_k,\, \widehat{\cV}_k) = \min_{\mQ \in \mathbb{O}(k)} \left( \max_{1 \leq j \leq n} \| \vv_j - \mQ\widehat{\vv}_j \|_2 \right), \label{eqn:twoinf_distance}
\end{equation}
where $\mathbb{O}(k)$ is the set of $k \times k$ orthogonal matrices. The minimization over $\mQ \in \mathbb{O}(k)$ makes the value of $d_\mathrm{row}$ independent of the particular choice of bases used to represent $\cV_k$ and $\widehat{\cV}_k$. Instead of measuring an aggregated total of subspace approximation errors, $d_\mathrm{row}$ measures how errors are distributed across the rows of two bases that have been optimally aligned with one another. \Cref{thm:idbound_randomized_normwise_componentwise} bounds the perturbations of $\varphi_\mathrm{max}$ induced by subspace errors measured under $d_\mathrm{row}$.
\begin{theorem} \label{thm:idbound_randomized_normwise_componentwise}
Choose $k$ such that $\sigma_k(\mA) > \sigma_{k + 1}(\mA)$. Given skeleton indices $J = \{j_1,\, \ldots,\, j_k \}$ and a $k$-dimensional subspace $\widehat{\cV}_k \subseteq \R^n$, let $\varphi_\mathrm{max}$ (resp. $\widehat{\varphi}_\mathrm{max}$) be the largest principal angle between $\cI_J$ and $\cV_k$ (resp. $\widehat{\cV}_k$). If $\widehat{\varphi}_\mathrm{max} < \pi / 2$, then
\begin{equation}
\cos \varphi_\mathrm{max} \geq \cos \widehat{\varphi}_\mathrm{max} - \frac{k c_k \mu}{\cos \widehat{\varphi}_\mathrm{max}} + \cO(\mu^2), \label{eqn:subspace_perturbation_componentwise}
\end{equation}
where $\mu = d_\mathrm{row}(\cV_k,\, \widehat{\cV}_k)$ and $c_k = \max_j \ell_j^{(k)}$ is the coherence of $\cV_k$.
\end{theorem}

\begin{proof}
Refer to \cref{subsec:idbound_randomized_componentwise_proof}.
\end{proof}

Note that the assumption $\widehat{\varphi}_\mathrm{max} < \pi / 2$ will always be satisfied by RGKS, due to the RRQR factorization in line \ref{line:rgks_rrqr} of the algorithm. For comparison with \cref{eqn:subspace_perturbation_componentwise}, the perturbation bound from \cref{thm:idbound_randomized_normwise} says that
\begin{equation}
\cos \varphi_\mathrm{max} \geq \cos \widehat{\varphi}_\mathrm{max} \cos \theta_\mathrm{max} - \sin \widehat{\varphi}_\mathrm{max} \sin \theta_\mathrm{max}. \label{eqn:subspace_perturbation_normwise}
\end{equation}
In a situation with large aggregate error but small component-wise error, we will have the approximate relationship $\mu \ll \sin \theta_\mathrm{max} \approx 1$. Then, provided that $k c_k$ is not too large and $\widehat{\varphi}_\mathrm{max}$ is not too close to $\pi / 2$, comparing equations \cref{eqn:subspace_perturbation_componentwise} and \cref{eqn:subspace_perturbation_normwise} shows that the component-wise errors measured by $\mu$ induce significantly smaller perturbations of $\cos \varphi_\mathrm{max}$ than an analysis using $\theta_\mathrm{max}$ would suggest.

To bound the value of $\mu$, we can draw on the work of authors who have previously developed subspace perturbation theory under errors measured by $d_\mathrm{row}$. For example, Damle and Sun \cite{damle_uniform_bounds} bound $d_\mathrm{row}$ between the eigenspaces of arbitrary symmetric matrices $\mM_1$ and $\mM_2$ which are symmetric perturbations of one another. Letting $\mM_1 = \mA\tp\mA$ and $\mM_2 = \widehat{\mA}\tp\widehat{\mA}$ results in a bound on $d_\mathrm{row}(\cV_k,\, \widehat{\cV}_k)$. The results of Zhang and Tang \cite{zhang_tang} provide $d_\mathrm{row}$ subspace perturbation bounds in the limit $n \to \infty$, under assumptions that $\| \mM_1 \|_2^{-1} \| \mM_1 - \mM_2 \|_2$ decays polynomially fast with the matrix dimension $n$. Cape et al.\ \cite{cape_two_infinity} provide further non-asymptotic bounds on $d_\mathrm{row}$. The recurring theme is that $d_\mathrm{row}$ tends to be smaller when $\cV_k$ is highly incoherent. 

\subsection{Empirical Analysis of RSVD Subspace Errors} \label{subsec:subspace_approximation_numerics}
This section investigates, empirically, how the RSVD approximation errors in line \ref{line:rgks_rsvd} of RGKS (\cref{alg:rgks}) are distributed across the components of $\cV_k$ and $\widehat{\cV}_k$, as well as how these errors affect the performance of RGKS. Whereas the previous section used $d_\mathrm{row}(\cV_k,\, \widehat{\cV}_k)$ to quantify component-wise errors, these experiments will focus on a related error measure: the element-wise discrepancies between $\mP$ and $\widehat{\mP}$, the orthogonal projection matrices onto $\cV_k$ and $\widehat{\cV}_k$. This allows for a more direct comparison with $\theta_\mathrm{max} = \arcsin \| \mP - \widehat{\mP} \|_2$, the ``aggregated'' subspace error measure that appears in \cref{thm:idbound_randomized_normwise}. We will show that even though $\theta_\mathrm{max}$ can approach $\pi/2$ when $\gamma_k \approx 1$, the component-wise errors $|(\mP - \widehat{\mP})_{ij}|$ often remain small, allowing for an accurate RGKS approximation in the end.

\begin{figure} [h]
    \centering
    \begin{subfigure}{.32\textwidth}
        \includegraphics[scale=.28]{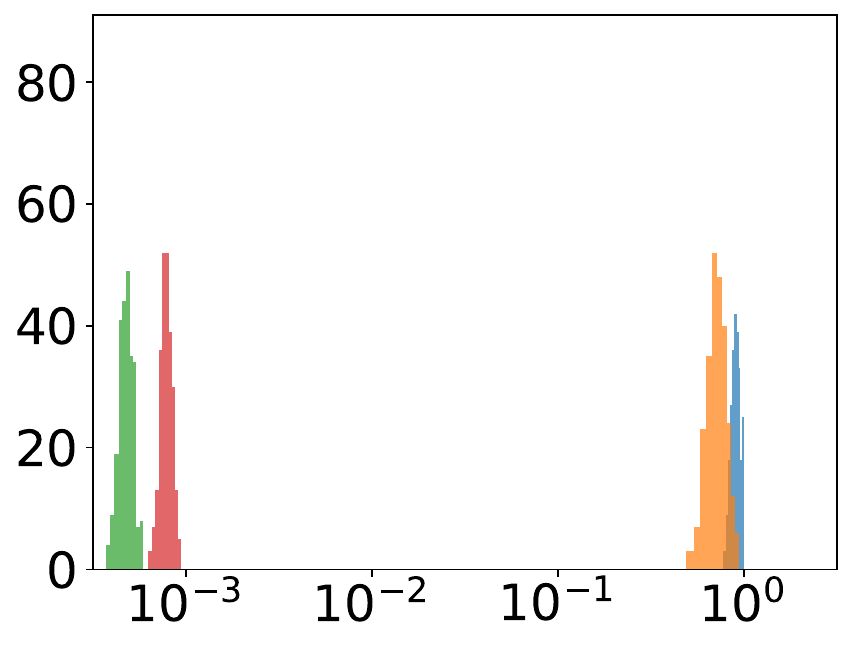}
        \caption{High coherence $\mV_k$.}
        \label{subfig:highcoherencehistogram}
    \end{subfigure}
    \hfill
    \begin{subfigure}{.32\textwidth}
        \includegraphics[scale=.28]{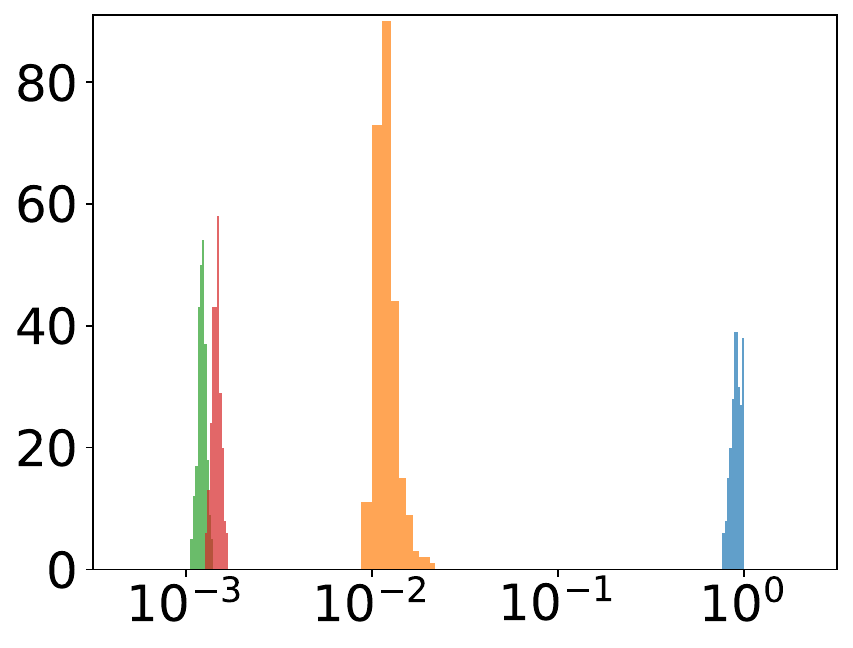}
        \caption{Low coherence $\mV_k$.}
        \label{subfig:lowcoherencehistogram}
    \end{subfigure}
    \hfill
    \begin{subfigure}{.32\textwidth}
        \includegraphics[scale=.28]{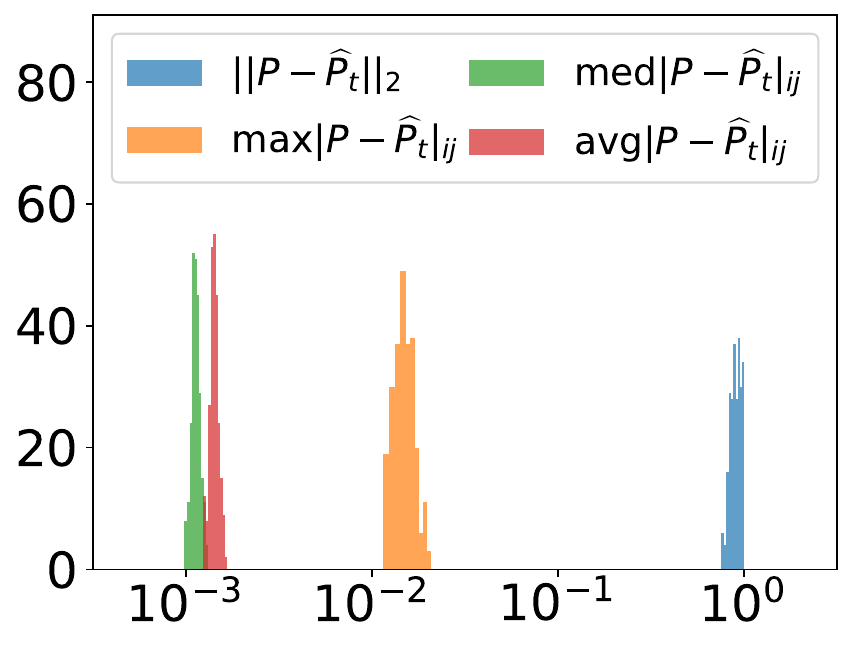}
        \caption{Random $\mV_k$.}
        \label{subfig:randomhistogram}
    \end{subfigure}
    \caption{Subspace approximation errors produced by RSVD and the maximum, median, and average element-wise error $|(\mP - \widehat{\mP})_{ij}|$ across $T=250$ trials. The experiment was run with approximation rank $k=70$ (for which $\gamma_k \approx 0.853 $ and $r_k\approx 254.020$), no power iteration ($q=0$), oversampling parameter $p=5$, and varying levels of coherence in $\mV_k$, the true right singular subspace of the $1024\times1024$ input matrix $\mA$.}
    \label{fig:metric_histograms}
\end{figure}

For these experiments, the spectral norm of $\mP-\widehat{\mP}$ as well as various statistics for the element-wise errors were recorded for matrices whose right singular subspaces  were derived from a noisy permutation matrix, a noisy dyadic Hadamard matrix, and a random Gaussian matrix. The approximation rank was chosen was chosen such that it lay at the tail end of a decay region of the singular spectrum. For ranks of this type, $\gamma_k$ is near unity, meaning that there will be significant noise in the RSVD subspace estimate. The results are shown in \cref{fig:metric_histograms}.
  
These experiments illustrate two distinct behaviors across difference subspace structure. For both low coherence and random right singular subspaces, the average, median, and maximum element-wise error are all orders of magnitude smaller than the subspace error as seen in \cref{subfig:lowcoherencehistogram} and \cref{subfig:randomhistogram}. Therefore, for matrices with with incoherent and random right singular subspaces we hypothesize that large $\theta_\mathrm{max}$ does not significantly affect the suitability of the approximation to be used in RGKS. In contrast, for coherent subspaces (see \cref{subfig:highcoherencehistogram}), the largest element-wise errors approach the same magnitude as the subspace errors. Nevertheless, the median and average element-wise errors are relatively small. Our hypothesis is that the largest of these element-wise errors occur along the diagonal of the projectors, corresponding to the leverage scores of $\cV_k$ and $\widehat{\cV}_k$. Note that for the highly coherent subspaces we constructed, the leverage score distribution of $\cV_k$ is mostly concentrated in a set of only $k$ indices, making it almost certain that these indices will be selected as skeleton columns even in the presence of noise. The supplementary text presents numerical experiments which demonstrate that even in the presence of noise, RGKS produces very high quality approximations of $\mA$ when the subspace is coherent.

\section{Conclusions}

We have showed how the accuracy of interpolative decomposition algorithms is affected by the properties of the input matrix being operated on, particularly with respect to singular value decay and singular subspace geometry. Motivated by these considerations, we introduced the RGKS algorithm, a novel interpolative decomposition which uses a randomized approximation of a singular subspace. Numerical experiments in \cref{section:numerics} showcased the myriad ways in which these interpolative decompositions are affected by singular value decay and singular subspace geometry, and simultaneously, demonstrated that RGKS was competitive with well-known algorithms for low-rank approximation. Surprisingly, these experiments showed that RGKS computed more accurate approximations than RSVD under certain circumstances. \Cref{section:analysis} presented error bounds which described the effects of input matrix structures on a generic interpolative decomposition, and \cref{section:rgks_analysis} specialized these error bounds to RGKS, as well as its deterministic counterpart GKS. In relating the GKS analysis to RGKS, we provided an analysis of randomization errors, while also finding that RGKS exhibits robustness to randomization noise which is difficult to explain with theory alone. We concluded with numerical experiments that shed light on the source of this robustness.

\section{Code Availability} Code to reproduce the results in this paper can be found at \url{https://github.com/robin-armstrong/structure-aware-id-analysis}.

\section{Acknowledgements}

RA and AD were partially supported by the National Science Foundation under award DMS-2146079. AD was also partially supported by the SciAI Center, funded by the Office of Naval Research under Grant Number N00014-23-1-2729. We would like to thank the two anonymous referees for their many helpful suggestions.

\bibliographystyle{siamplain}
\bibliography{sources}

\appendix

\section{Proofs} \label{section:proofs}

\subsection{Proof of Conditioning-Based Error Bound} \label{subsec:lowrankbound_spectrumonly_proof}

Let $\cW \subseteq \R^m$ be a subspace with $\dim \cW < m$, and let $\mP_\cW$ be the orthogonal projector onto $\cW$. \Cref{thm:lowrankbound_spectrumonly} states that if $\mE = \mA - \mP_\cW\mA$ and $\sigma_{k + 1}(\mE) > 0$, then $\| \mE \|_2 \leq \| \mSigma_\perp \|_2 \kappa(\mE,\, k + 1)$, where $\kappa(\mE,\, k + 1) = \sigma_1(\mE) / \sigma_{k + 1}(\mE)$ is a modified condition number. The proof uses an Ostrowsky-type singular value bound \cite[Theorem 6.1]{ipsen_1998}, which lets us write $\sigma_{k + 1}(\mE) = \sigma_{k + 1}((\mI - \mP_\cW)\mA) \leq  \| \mSigma_\perp \|_2 \| \mI - \mP_\cW \|_2$. Since $\mI - \mP_\cW$ is the orthogonal projector onto $\cW^\perp$, the assumption that $\dim \cW < m$ implies $\| \mI - \mP_\cW \|_2 = 1$. Thus $\sigma_{k + 1}(\mE) \leq \| \mSigma_\perp \|_2$, and since $\| \mSigma_\perp \|_2 = \sigma_{k + 1}(\mA) > 0$, this implies that
\begin{equation*}
\| \mE \|_2 = \| \mSigma_\perp \|_2 \frac{\sigma_1(\mE)}{\| \mSigma_\perp \|_2} \leq \| \mSigma_\perp \|_2 \frac{\sigma_1(\mE)}{\sigma_{k + 1}(\mE)} = \| \mSigma_\perp \|_2 \kappa(\mE,\, k + 1).
\end{equation*}

\subsection{Proof of Principal Angle Lemma} \label{subsec:angle_svd_connection_proof}

To allow the use of the CS decomposition, we impose the restriction that $k \leq n/2$. Given a selection $J = \{ j_1,\, \ldots,\, j_k \}$ of column indices, let $\varphi_1,\, \ldots,\, \varphi_k$ be the principal angles between $\cI_J = \range(\mI_{:,\, J})$ and $\cV_k = \range(\mV_k)$. The first claim of \cref{lemma:angle_leverage_connection} is that $\cos \varphi_1,\, \ldots,\, \cos \varphi_k$ are the singular values of $\mV_{J,\, 1:k}$. This follows from the fact that $\mI_{:,\, J}$ and $\mV_k$ are orthonormal bases for $\cI_J$ and $\cV_k$, and that $\mV_{J,\, 1:k} = (\mI_{:,\, J})\tp \mV_k$. We now have that $\sigma_\mathrm{min}(\mV_{J,\, 1:k}) = \cos \varphi_\mathrm{max}$, which implies that $\mV_{J,\, 1:k}$ is invertible if and only if $\max_i \varphi_i < \pi / 2$. This proves the second claim of \cref{lemma:angle_leverage_connection}.   

The final claim of \cref{lemma:angle_leverage_connection} is that the singular values of $\mV_{[n] \setminus J,\, 1:k}(\mV_{J,\, 1:k})^{-1}$ are equal to $\tan \varphi_1,\, \ldots,\, \tan \varphi_k$, provided $\mV_{J,\, 1:k}$ is invertible. Let $\mPi$ be a permutation which moves indices $J$ to the front (i.e., $\mPi = \begin{bmatrix} \mI_{:,\, J} & \mI_{:,\, [n] \setminus J} \end{bmatrix}$), let $\mV = [\mV_k \:\: \mV_\perp]$, and consider a CS decomposition of $\mPi\tp\mV$:
\begin{equation*}
\begin{bmatrix}
\mV_{J,\, 1:k} & \mV_{J,\, k + 1:n} \\
\mV_{[n] \setminus J,\, 1:k} & \mV_{[n] \setminus J,\, k + 1:n}
\end{bmatrix} = \begin{bmatrix}
\mQ_1 \mC \mW_1\tp & \mQ_1 \begin{bmatrix} \mS & \mZero \end{bmatrix} \mW_2\tp \\
\mQ_2 \begin{bmatrix} \mS \\ \mZero \end{bmatrix} \mW_1\tp & \mQ_2 \begin{bmatrix} -\mC & \mZero \\ \mZero & \mI \end{bmatrix} \mW_2\tp
\end{bmatrix},
\end{equation*}
where $\mQ_1,\, \mW_1 \in \R^{k \times k}$ and $\mQ_2,\, \mW_2 \in \R^{(n - k) \times (n - k)}$ are orthogonal matrices, $\mC = \mathrm{diag}(\cos \varphi_1,\, \ldots,\, \cos \varphi_k)$, and $\mS = \mathrm{diag}(\sin \varphi_1,\, \ldots,\, \sin \varphi_k)$. The claim follows from
\begin{align}
\mV_{[n] \setminus J,\, 1:k}(\mV_{J,\, 1:k})^{-1} &= \mQ_2 \begin{bmatrix} \mS\mC^{-1} \\ \mZero \end{bmatrix} \mQ_1\tp \nonumber \\
(\mV_{J,\, 1:k})^{-1}\mV_{J,\, k + 1:n} &= \mW_1 \begin{bmatrix} \mC^{-1}\mS & \mZero \end{bmatrix} \mW_2\tp, \label{eqn:cs_tangent_equation}
\end{align}

where $\mS\mC^{-1} = \mC^{-1}\mS = \mathrm{diag}(\tan \varphi_1,\, \ldots,\, \tan \varphi_k)$.

\subsection{Proof of Spectral Norm Subspace Geometry Bound}
\label{subsec:idbound_subspace_proof}

\Cref{thm:idbound_subspaceonly} states that if $\varphi_\mathrm{max} \defeq \max_i \varphi_i < \pi / 2$, then $\| \mA - \mA_{:,\, J}(\mA_{:,\, J})\pinv\mA \|_2 \leq \| \mSigma_\perp \|_2 \sec \varphi_\mathrm{max}$. 
To prove this inequality, partition $\mPi\tp\mV$ as
\begin{equation}
\mPi\tp\mV = \begin{bNiceMatrix}[last-row, last-col]
\mV_{11} & \mV_{12} & \mbox{\scriptsize $k$} \\
\mV_{21} & \mV_{22} & \mbox{\scriptsize $n - k$} \\
\mbox{\scriptsize $k$} & \mbox{\scriptsize $n - k$}
\end{bNiceMatrix}\, . \label{eqn:V_partition}
\end{equation}
Let $\| \cdot \|$ denote the spectral or Frobenius norm. By applying \cite[Theorem 9.1]{halko_finding_structure_with_randomness} with $\mOmega = \mPi_{:,\, 1:k},\, \mOmega_1 = \mV_{11}\tp$, and $\mOmega_2 = \mV_{12}\tp$, we find that
\begin{equation}
\| \mA - \mA_{:,\, J}(\mA_{:,\, J})\pinv\mA \|^2 \leq \| \mSigma_\perp \|^2 + \| \mSigma_\perp \|_2^2 \| \mV_{11}^{-1}\mV_{12} \|^2,
\label{eqn:hmt2011_idversion}
\end{equation}
provided that $\mV_{11}$ is full rank. In light of \cref{lemma:angle_leverage_connection}, this is equivalent to requiring that $\varphi_\mathrm{max} < \pi / 2$. Equation \cref{eqn:cs_tangent_equation} shows that the singular values of $\mV_{11}^{-1}\mV_{12} = (\mV_{J,\, 1:k})^{-1}\mV_{J,\, k + 1:n}$ are equal to $\tan \varphi_1,\, \ldots,\, \tan \varphi_k$. The spectral norm of $\mV_{11}^{-1}\mV_{12}$ is therefore $\tan \varphi_\mathrm{max}$, meaning that
\begin{equation*}
\| \mA - \mA_{:,\, J}(\mA_{:,\, J})\pinv\mA \|_2^2 \leq \| \mSigma_\perp \|_2^2(1 + \tan^2 \varphi_\mathrm{max}) = \| \mSigma_\perp \|_2^2 \sec^2 \varphi_\mathrm{max}.
\end{equation*}

\subsection{Proof of Frobenius Norm Subspace Geometry Bound} \label{subsec:idbound_bothstructures_proof}

\Cref{thm:idbound_bothstructures} states that 
\[\| \mA - \mA_{:,\, J}(\mA_{:,\, J})\pinv\mA \|_\frob \leq \| \mSigma_\perp \|_\frob \sqrt{1 + \frac{1}{r_k} \sum_{i = 1}^k \tan^2 \varphi_i},
\]
where $r_k = \| \mSigma_\perp \|_\frob^2 / \| \mSigma_\perp \|_2^2$ is the residual stable rank, provided that $\varphi_\mathrm{max} < \pi / 2$. The proof begins with \cref{eqn:hmt2011_idversion}, taking $\| \cdot \|$ to be the Frobenius norm. We then insert the definition of $r_k$ and use \cref{eqn:cs_tangent_equation} to expand $\| \mV_{11}^{-1}\mV_{12} \|_\frob$ in terms of $\tan \varphi_1,\, \ldots,\, \tan \varphi_k$.

\subsection{Proof of Full Subspace Perturbation Bound} \label{subsec:idbound_randomized_proof}

The statement in \cref{thm:idbound_randomized_normwise} is equivalent to $\cos \varphi_\mathrm{max} \geq \cos(\widehat{\varphi}_\mathrm{max} + \theta_\mathrm{max})$. As discussed in the proof of \cref{lemma:angle_leverage_connection}, $\cos \varphi_\mathrm{max} = \sigma_\mathrm{min}(\mV_{J,\, 1:k}) = \sigma_\mathrm{min}(\mV_{11})$, following the partition of $\mPi\tp\mV$ given in \cref{eqn:V_partition}. Similarly, letting $\widehat{\mV}_k \in \R^{n \times k}$ and $\widehat{\mV}_\perp \in \R^{n \times (n - k)}$ be orthonormal bases for $\widehat{\cV}_k$ and $\widehat{\cV}_k^\perp$, we define $\widehat{\mV} = [ \widehat{\mV}_k \:\: \widehat{\mV}_\perp ]$. Consider the CS decomposition
\begin{equation*}
\widehat{\mV}\tp\mV = \begin{bmatrix}
\mQ_1 \mC \mW_1\tp & \mQ_1 \begin{bmatrix} \mS & \mZero \end{bmatrix} \mW_2\tp \\
\mQ_2 \begin{bmatrix} \mS \\ \mZero \end{bmatrix} \mW_1\tp & \mQ_2 \begin{bmatrix} -\mC & \mZero \\ \mZero & \mI \end{bmatrix} \mW_2\tp
\end{bmatrix},
\end{equation*}
where $\mQ_1,\, \mW_1 \in \R^{k \times k}$ and $\mQ_2,\, \mW_2 \in \R^{(n - k) \times (n - k)}$ are orthogonal matrices, $\mC = \mathrm{diag}(\cos \theta_1,\, \ldots,\, \cos \theta_k)$, $\mS = \mathrm{diag}(\sin \theta_1,\, \ldots,\, \sin \theta_k)$, and $\theta_1,\, \ldots,\, \theta_k$ are the principal angles between $\cV_k$ and $\widehat{\cV}_k$. Then 
\[
\mV_k = \widehat{\mV}(\widehat{\mV}\tp\mV)_{:,\, 1:k} = \widehat{\mV}_k \mQ_1\mC\mW_1\tp + \widehat{\mV}_\perp\mQ_2\mPsi\mW_1\tp,
\]
where $\mPsi = [\mS \:\: \mZero]\tp$. Define $\widehat{\mV}_{ij}$ for $i,\, j \in \{ 1,\, 2 \}$ by partitioning $\mPi\tp\widehat{\mV}$ in a way identical to \cref{eqn:V_partition}, so that $\cos \widehat{\varphi}_\mathrm{max} = \sigma_\mathrm{min}(\widehat{\mV}_{11})$ and $\sin \widehat{\varphi}_\mathrm{max} = \sigma_\mathrm{max}(\widehat{\mV}_{12})$. We now have,
\begin{align*}
\mV_{11} = (\mPi_{:,\, 1:k})\tp\mV_k &= (\mPi_{:,\, 1:k})\tp \widehat{\mV}_k \mQ_1\mC \mW_1\tp + (\mPi_{:,\, 1:k})\tp \widehat{\mV}_\perp\mQ_2\mPsi \mW_1\tp \\
&= (\widehat{\mV}_{11}\mQ_1\mC + \widehat{\mV}_{12}\mQ_2\mPsi)\mW_1\tp.
\end{align*}
Using the orthogonal invariance of singular values, an Ostrowsky-type bound \cite[Theorem 6.1]{ipsen_1998}, and Weyl's inequality,
\begin{align*}
\cos \varphi_\mathrm{max} = \sigma_\mathrm{min}(\mV_{11}) &= \sigma_\mathrm{min}(\widehat{\mV}_{11}\mQ_1\mC + \widehat{\mV}_{12}\mQ_2\mPsi) \\
&\geq \sigma_\mathrm{min}(\widehat{\mV}_{11})\cos \theta_\mathrm{max} - \sigma_\mathrm{max}(\widehat{\mV}_{12}) \sin \theta_\mathrm{max} \\
&= \cos \widehat{\varphi}_\mathrm{max} \cos \theta_\mathrm{max} - \sin \widehat{\varphi}_\mathrm{max} \sin \theta_\mathrm{max} \\
&= \cos(\widehat{\varphi}_\mathrm{max} + \theta_\mathrm{max}),
\end{align*}
which completes the proof.

\subsection{Proof of Restricted Subspace Geometry Error Bound} \label{subsec:rgks_smallerskeleton_proof}

Let $\widehat{\cV}_{k - t}$ be the subspace spanned by the first $k - t$ right singular vectors estimated in line \ref{line:rgks_rsvd} of RGKS (\cref{alg:rgks}), let $\widehat{\varphi}_\mathrm{max}^{(I)}$ be the largest principal angle between $\cI_I$ and $\widehat{\cV}_{k-t}$ for some $I \subseteq J$, and let $\theta_\mathrm{max}\ps{k - t}$ be the largest principal angle between $\cV_{k - t}$ and $\widehat{\cV}_{k - t}$. Equation \cref{eqn:rgksbound_smallerskeleton} states that
\begin{equation*}
\| \mA - \mA_{:,\, J}(\mA_{:,\, J})\pinv\mA \|_2 \leq \frac{\| \mSigma_\perp \|_2 \sec(\widehat{\varphi}_\mathrm{max}^{(I)} + \theta_\mathrm{max}\ps{k - t})}{\gamma_{k - t + 1,\, k + 1}}.
\end{equation*}
To prove this bound, consider an arbitrary subset $I \subseteq J$ of size $k - t$. If $\widehat{\varphi}_\mathrm{max}\ps{I}$ is the largest angle between $\widehat{\cV}_{k - t}$ and $\cI_I$, and $\theta_\mathrm{max}\ps{k - t}$ is the largest angle between $\cV_{k - t}$ and $\widehat{\cV}_{k - t}$, then we have
\begin{equation*}
\varphi_\mathrm{max}\ps{I} \leq \widehat{\varphi}_\mathrm{max}\ps{I} + \theta_\mathrm{max}\ps{k - t}
\end{equation*}
This is a direct consequence of \cref{thm:idbound_randomized_normwise}, and inserting this into \cref{eqn:idbound_subspaceonly_subsetangle} completes the proof.

\subsection{Proof of Row-Wise Subspace Perturbation Bound} \label{subsec:idbound_randomized_componentwise_proof}

We now let $\vv_j$ and $\widehat{\vv}_j$ denote the $j\nth$ rows of $\mV_k$ and $\widehat{\mV}_k$, respectively, for $j = 1,\, \ldots,\, n$. \Cref{thm:idbound_randomized_normwise_componentwise} states that if $\mu = d_\mathrm{row}(\cV_k,\, \widehat{\cV}_k)$ and $\widehat{\varphi}_\mathrm{max} < \pi / 2$, then
\begin{equation*}
\cos \varphi_\mathrm{max} \geq \cos \widehat{\varphi}_\mathrm{max} - \frac{k c_k \mu}{\cos \widehat{\varphi}_\mathrm{max}} + \cO(\mu^2),
\end{equation*}
where $c_k = \max_j \ell_j^{(k)}$ is the coherence of $\cV_k$. To prove this, let $\mV_{11} = \mV_{J,\, 1:k}$ and $\widehat{\mV}_{11} = \widehat{\mV}_{J,\, 1:k}$, as in the proof of \cref{thm:idbound_randomized_normwise}. Using Weyl's inequality,
\begin{equation}
\sigma_k(\mV_{11}\mV_{11}\tp) \geq \sigma_k(\widehat{\mV}_{11}\widehat{\mV}_{11}\tp) - \| \mV_{11}\mV_{11}\tp - \widehat{\mV}_{11}\widehat{\mV}_{11}\tp \|_\frob. \label{eqn:componentbound_weyl_frobenius}
\end{equation}
Each element of $\mV_{11}\mV_{11}\tp - \widehat{\mV}_{11}\widehat{\mV}_{11}\tp$ has the form $\vv_i\tp\vv_j - \widehat{\vv}_i\tp\widehat{\vv}_j$ for some $i,\, j \in J$. These elements can be bounded in magnitude as follows: let $\mQ_0$ be a $k \times k$ orthogonal matrix achieving the minimum in \cref{eqn:twoinf_distance}. We then have $\| \vv_s - \mQ_0\widehat{\vv}_s \|_2 \leq \mu$ for $1 \leq s \leq n$, and along with $\| 
\vv_s \|_2 = \ell_s^{(k)} \leq c_k$, this implies that
\begin{align*}
|\vv_i&\tp\vv_j - \widehat{\vv}_i\tp\widehat{\vv}_j| \\
&= |\vv_i\tp\vv_j - (\mQ_0\widehat{\vv}_i)\tp(\mQ_0\widehat{\vv}_j)| \\
&= |\vv_i\tp(\mQ_0 \widehat{\vv}_j - \vv_j) + \vv_j\tp(\mQ_0\widehat{\vv}_i - \vv_i) + (\mQ_0\widehat{\vv}_i - \vv_i)\tp(\mQ_0\widehat{\vv}_j - \vv_j)| \\
&\leq \| \vv_i \|_2 \| \vv_j - \mQ_0\widehat{\vv}_j \|_2 + \| \vv_j \|_2\| \vv_i - \mQ_0\widehat{\vv}_i \|_2 + \| \vv_i - \mQ_0\widehat{\vv}_i \|_2\| \vv_j - \mQ_0\widehat{\vv}_j \|_2 \\
&\leq 2c_k\mu + \mu^2.
\end{align*}
From this bound, it follows that $\| \mV_{11}\mV_{11}\tp - \widehat{\mV}_{11}\widehat{\mV}_{11}\tp \|_\frob \leq 2kc_k\mu + k\mu^2$.

Returning to equation \cref{eqn:componentbound_weyl_frobenius}, \cref{lemma:angle_leverage_connection} shows that $\cos \varphi_\mathrm{max} = \sigma_k(\mV_{11})$, and similarly, $\cos \widehat{\varphi}_\mathrm{max} = \sigma_k(\widehat{\mV}_{11})$. This implies $\sigma_k(\mV_{11}\mV_{11}\tp) = \cos^2 \varphi_\mathrm{max}$ and $\sigma_k(\widehat{\mV}_{11}\widehat{\mV}_{11}\tp) = \cos^2 \widehat{\varphi}_\mathrm{max}$, and therefore, $\cos^2 \varphi_\mathrm{max} \geq \cos^2 \widehat{\varphi}_\mathrm{max} - 2kc_k\mu - k\mu^2$. Thus $\cos \varphi_\mathrm{max} \geq f(\eta)$, where $\eta = 2kc_k\mu + k\mu^2$ and $f(x) = \sqrt{\cos^2 \widehat{\varphi}_\mathrm{max} - x}$. We conclude by writing $f(\eta) = f(0) + \eta f'(0) + \cO(\eta^2)$ and re-inserting the definition of $\eta$.

\end{document}


\maketitle

This document offers figure and analyses which supplement the main text. As in the main text, $\mA$ denotes an $m \times n$ real matrix, and $k$ denotes the target rank of a low-rank approximation of $\mA$. We write $\mA = \mU_k\mSigma_k\mV_k\tp + \mU_\perp\mSigma_\perp\mV_\perp\tp$ to denote a singular value decomposition partitioned at rank $k$, so that $\mU_k\mSigma_k\mV_k\tp$ is the optimal rank-$k$ approximation of $\mA$ in the spectral or Frobenius norm, and $\| \mSigma_\perp \|$ is the corresponding approximation error. We write $\mA \approx \widehat{\mU}_k\widehat{\mSigma}_k\widehat{\mV}_k\tp$ to denote a randomized SVD of $\mA$ with $k$ components.

\section{Computational Efficiency of RGKS}

\begin{figure}[H]
    \centering
    \includegraphics[scale=.52]{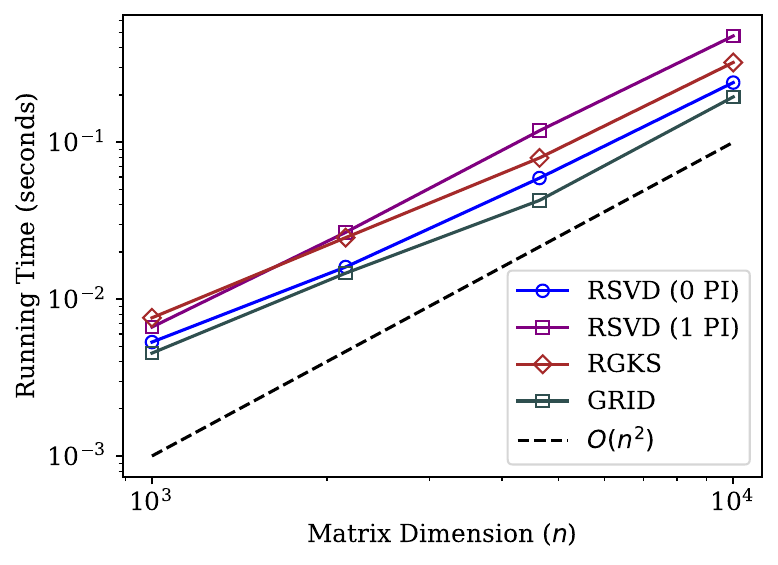}
    \caption{Runtime comparison between low-rank approximation algorithms applied to four square matrices with dimensions ranging from 1,000 to 10,000. Each data point is the average runtime over 10 low-rank approximations of the same matrix. Tests were run using Julia on a 2.9 GHz Intel i7-10700 CPU, with 16 GB of available RAM. All algorithms were precompiled with a dry-run before making runtime measurements.}
    \label{fig:runtime_comparison}
\end{figure}

RGKS is more efficient than several algorithms that are similar to it in design. Notably, it is more efficient than the standard approach of selecting columns via a RRQR factorization of $\mA$, since an RRQR factorization of $\mA$ operates on columns of length $m$, whereas RGKS considers only the rows of $\widehat{\mV}_k$ with length $k \ll m$. This is similar to the RID algorithm \cite{martinsson_rokhlin_tygert_randid, woolfe_srft}, \add{referred to as RCPQR in the main text}, which selects skeleton columns using by applying an RRQR factorization to a \delete{Gaussian-}random linear embedding of the columns of $\mA$. RGKS can be expected to have slightly longer runtimes than RID, since the former uses matrix-vector products with both $\mA$ and $\mA\tp$, while the latter uses only products products with $\mA\tp$. \Cref{fig:runtime_comparison} plots the runtime of RGKS in comparison to RSVD with and without power iteration and RID. It shows that in terms of wall-clock time, all these algorithms are similar to one another, and that when applied to $n \times n$ matrices, they have the same asymptotic complexity $\cO(n^2)$ for fixed $k$.

\section{Numerical Comparison of Spectral Norm Errors}

This section presents spectral norm counterparts to the experiments in section 5 of the main text. 

\begin{figure}[H]
    \centering
    \begin{subfigure}{.95\textwidth}
        \includegraphics[scale=.55]{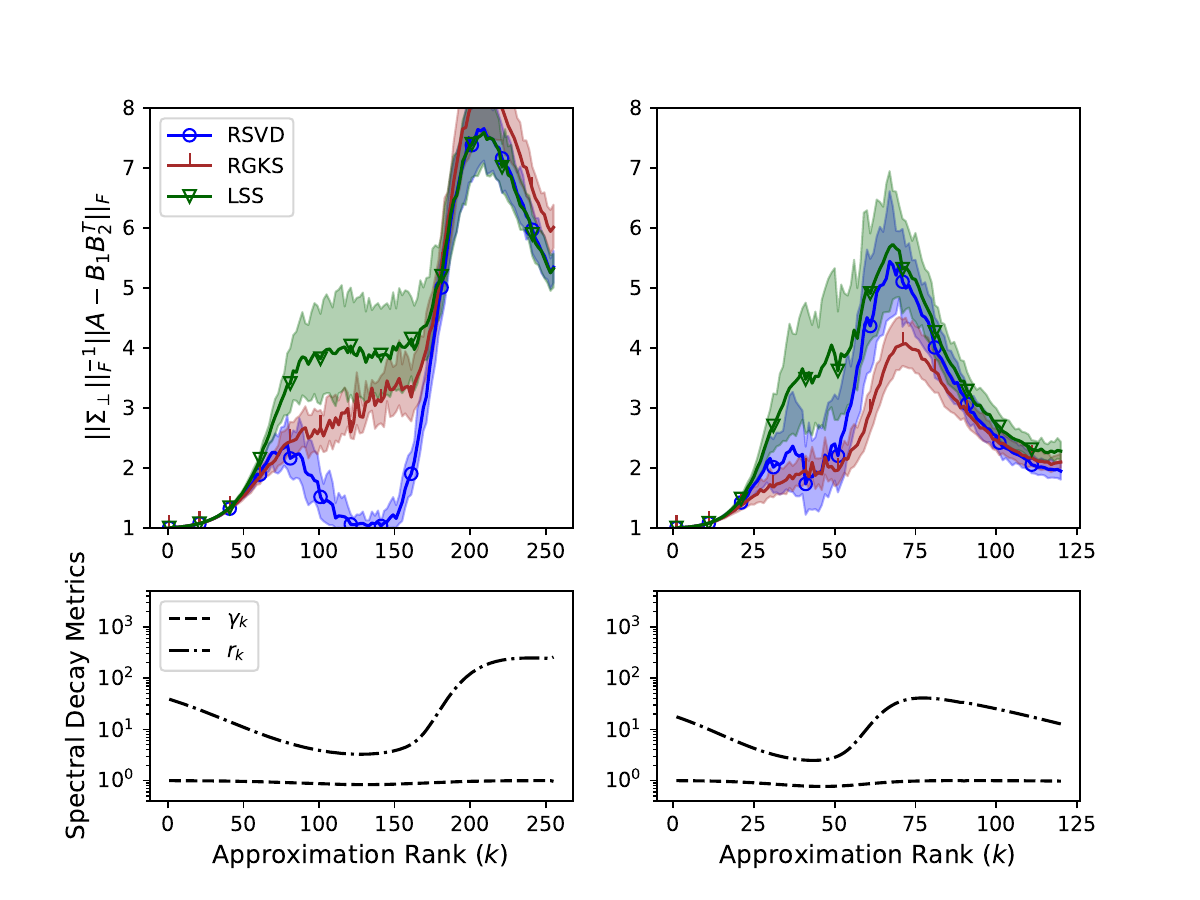}
        \caption{Spectral norm errors for  $256 \times 256$ test matrices with identical right singular subspaces.}
        \label{fig:supplemental_spectraleffects}
    \end{subfigure}
    \hfill
    \begin{subfigure}{.95\textwidth}
        \includegraphics[scale=.55]{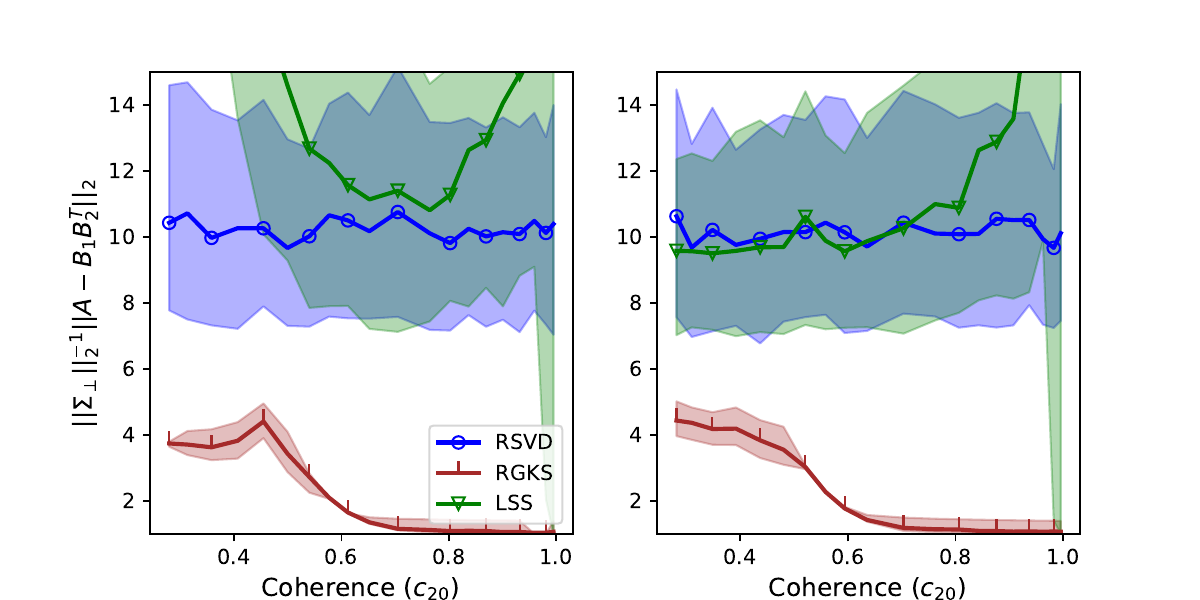}
        \caption{Spectral norm errors for square matrices of dimension 256 (left) and 252 (right) at rank $k = 20$. All of the matrices' singular spectra were identical in the first 252 entries.}
        \label{fig:supplemental_subspaceeffects}
    \end{subfigure}
    \caption{Spectral norm comparison of algorithm errors. Each data point is the mean error over 100 approximations of the same matrix. \add{Note that LSS refers to the leverage score sampling algorithm, also called LEVG in the main text.}}
\end{figure}

\section{RGKS Approximation Errors for Coherent Subspaces}
\begin{figure}[H]
    \centering
    \includegraphics[scale=.35]{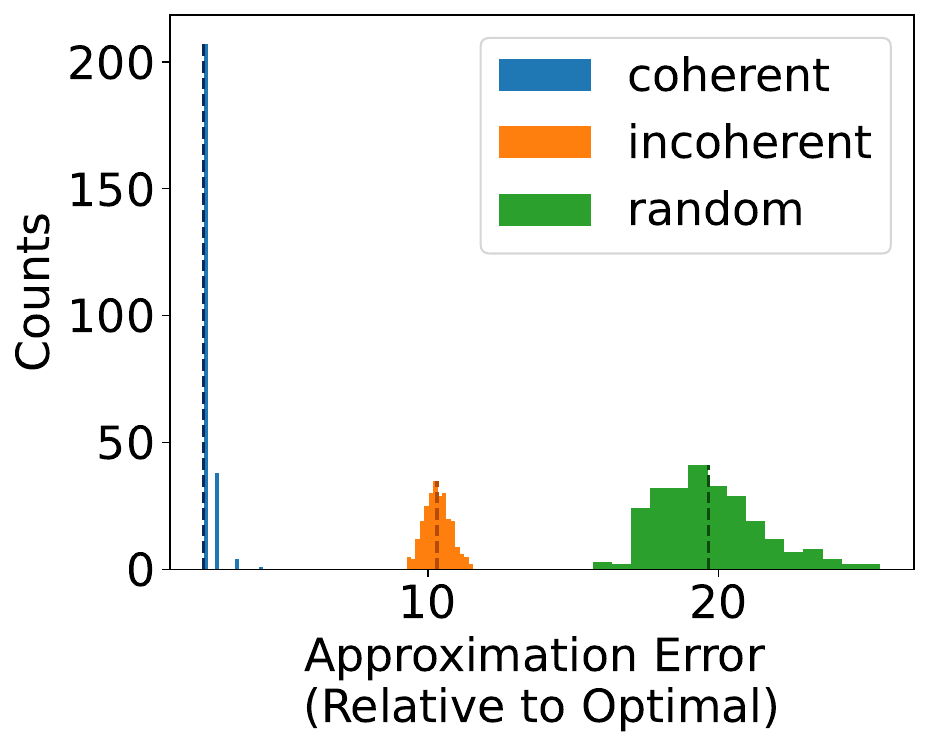}
        \caption{Histogram of RGKS approximation error relative to optimal. For right singular subspaces of the input matrix with varying levels of coherence, rank $k=70$ approximations were computed using oversampling $p=5$ and no power iteraton ($q=0$) across $T=250$ trials. The dashed lines represent the approximation error relative to the optimal when the deterministic GKS algorithm is used.}
        \label{fig:approximationerrors}
\end{figure}

Section 7.4 in the main text demonstrates, through numerical experiments, that certain element-wise differences between the true projector $\mP = \mV_k\mV_k\tp$ and the RSVD-approximated projector $\widehat{\mP} = \widehat{\mV}_k\widehat{\mV}_k\tp$ can be large when $\mV_k$ is highly coherent, and that all element-wise differences tend to be smaller when $\mV_k$ is less coherent. \Cref{fig:approximationerrors} demonstrates how these differences affect the approximation error of RGKS for coherent, incoherent, and randomly generated\footnote{This refers to instances of $\mV_k$ obtained by orthogonalizing a Gaussian random matrix.} subspaces. The approximation error is near the optimal value as well as the error that the deterministic GKS algorithm produces, despite the large element-wise errors present in the coherent case. In the incoherent and random cases, the RGKS error is distributed around the deterministic GKS error value, demonstrating that the error of the randomized version of the algorithm are not entirely one-sided.

\section{Near-Optimality of Condition-Number Bound}

Consider a subspace $\cW \subseteq \R^m$ with $\dim \cW = k < m$, defining a rank-$k$ approximation $\mA \approx \mP_\cW\mA$, where $\mP_\cW$ is the orthogonal projection matrix onto $\cW$. Letting $\mE = \mA - \mP_\cW\mA$, Theorem 6.4 in the main text states that
\begin{equation}
\| \mE \|_2 \leq \| \mSigma_\perp \|_2 \kappa(\mE,\, k + 1), \label{eqn:condition_number_tightness}
\end{equation}
where $\kappa(\mE,\, k + 1) = \sigma_1(\mE) / \sigma_{k + 1}(\mE)$ is a modified condition number. Figure 6.3 in the main text shows that for certain values of $k$ that correspond to a small or non-existent singular value gap, this bound is strikingly tight. This is due to the inequality
\begin{equation*}
\sigma_i(\mE) \geq \sigma_{k + i}(\mA),
\end{equation*}
which holds for $1 \leq i \leq \min\{ m,\, n \} - k$. To show this inequality, write $\mE = \mE_{i - 1} + \mF_{i - 1}$, where $\mE_{i - 1}$ is an optimal rank-$(i - 1)$ approximation of $\mE$ in the sense of the Eckart-Young theorem, and $\mF_{i - 1}$ is the residual of this approximation. Then $\sigma_i(\mE) = \| \mF_{i - 1} \|_2$. Furthermore, the relations $\rank \mP_\cW\mA \leq k$ and $\rank \mE_{i - 1} \leq i - 1$ imply that
\begin{equation}
\rank(\mP_\cW\mA + \mE_{i - 1}) \leq k + i - 1. \label{eqn:projected_rank_bound}
\end{equation}
By Weyl's inequality,
\begin{align*}
\sigma_{k + i}(\mA) &= \sigma_{k + i}(\mP_\cW\mA + \mE_{i - 1} + \mF_{i - 1}) \\
&\leq \sigma_{k + i}(\mP_\cW\mA + \mE_{i - 1}) + \| \mF_{i - 1} \|_2 \\
&= \sigma_{k + i}(\mP_\cW\mA + \mE_{i - 1}) + \sigma_i(\mE).
\end{align*}
Equation \cref{eqn:projected_rank_bound} shows that $\sigma_{k + i}(\mP_\cW \mA + \mE_{i - 1}) = 0$, completing the proof.

Suppose, now, that the singular spectrum of $\mA$ has very slow decay (or no decay), so that $\sigma_{k + 1}(\mA) \leq (1 + \varepsilon)\sigma_{2k + 1}(\mA)$ for some small number $\varepsilon > 0$. Assuming that $2k + 1 \leq \min\{m,\, n\}$, \cref{eqn:condition_number_tightness} shows that $\sigma_{k + 1}(\mE) \geq \sigma_{2k + 1}(\mA)$. We then have
\begin{equation*}
\| \mE \|_2 = \| \mSigma_\perp \|_2 \frac{\sigma_1(\mE)}{\sigma_{k + 1}(\mA)} \geq \| 
\mSigma_\perp \|_2 \frac{\sigma_1(\mE)}{(1 + \varepsilon)\sigma_{2k + 1}(\mA)} \geq \| 
\mSigma_\perp \|_2 \frac{\kappa(\mE,\, k + 1)}{1 + \varepsilon},
\end{equation*}
meaning that the bound is tight to within a factor of $1 + \varepsilon$.

\bibliographystyle{siamplain}
\bibliography{sources}